\documentclass[11pt,reqno]{amsart}

\usepackage[T1]{fontenc}
\usepackage[margin=1in]{geometry}
\usepackage{amsmath,amssymb,mathtools,mathrsfs}
\usepackage{microtype}
\usepackage{needspace}
\usepackage[hidelinks]{hyperref}
\hypersetup{
	pdftitle={Second class particle in ASEP},
	pdfauthor={Alexey Bufetov and Adyan Dordzhiev},
	pdfsubject={Second class particle in ASEP}
}

\numberwithin{equation}{section}
\newtheorem{theorem}{Theorem}[section]
\newtheorem{proposition}[theorem]{Proposition}
\newtheorem{lemma}[theorem]{Lemma}

\theoremstyle{definition}
\newtheorem{definition}[theorem]{Definition}
\newtheorem{example}[theorem]{Example}
\theoremstyle{remark}
\newtheorem{remark}[theorem]{Remark}

\usepackage[dvipsnames]{xcolor}
\usepackage[normalem]{ulem}

\newcommand{\PP}{\mathbb P}
\newcommand{\EE}{\mathbb E}
\newcommand{\ZZ}{\mathbb Z}
\newcommand{\NN}{\mathbb Z_{\ge0}}
\newcommand{\ind}{\mathbf 1}
\newcommand{\HH}{\mathcal H_q}
\newcommand{\T}{\mathbf T}
\newcommand{\id}{\mathrm{id}}
\newcommand{\up}{\mathcal U}
\newcommand{\soft}{\mathsf Z_q}

\newcommand{\distq}{\mathsf d_q}
\newcommand{\comp}{\mathsf C}
\newcommand{\poly}{\mathsf P}

\allowdisplaybreaks[1]
\title[Second-class particle in ASEP]{Second-class particle in ASEP. Rarefaction fan}
\author{Alexey Bufetov}
\address[Alexey Bufetov]{Institute of Mathematics,
Leipzig University, Germany}
\email{alexey.bufetov@gmail.com}

\author{Adyan Dordzhiev}
\address[Adyan Dordzhiev]{Max Planck Institute for Mathematics
in the Sciences, Leipzig, Germany}
\email{aedordzhiev@gmail.com}
\date{}
\begin{document}
	
	\begin{abstract}	
We study the asymptotic behavior of a second class particle in ASEP on the integer lattice. We prove a general formula for the distribution of its position at any fixed time for arbitrary initial conditions. As a first application, we obtain its limiting distribution for arbitrary initial conditions in the rarefaction fan case, generalizing the result of Cator--Pimentel from TASEP to ASEP.

	\end{abstract}
	\maketitle
	\tableofcontents
	
	\section{Introduction}\label{sec:introduction}
	\subsection{The model}\label{sec:intro-second-class}
The asymmetric simple exclusion process (ASEP) is one of the basic
one-dimensional interacting particle systems; see, e.g.,
Liggett~\cite{Liggett1999} for an introduction to the model.
It is a continuous-time Markov process
$(\eta_t)_{t\geq 0}$ on $\{0,1\}^{\mathbb{Z}}$,
where $\eta_t(x)=1$ means that site $x$ is occupied by a particle and $\eta_t(x)=0$ means that it is empty. Each particle independently attempts to jump one step
to the right at rate $1$ and one step to the left at rate $q$, where we assume $0<q<1$. A jump is suppressed if the target site is occupied. The process can be constructed by placing two independent
Poisson clocks at each site, with rates $1$ and $q$ for right and left jumps, respectively. We write
\begin{equation}\label{eq:intro-gamma-flux}
    \gamma=1-q,\qquad j(\rho)=\gamma\rho(1-\rho)
\end{equation}
for the asymmetry parameter and the hydrodynamic flux
of the system.

Our main object of study is a second-class particle. It can be introduced through the basic coupling, in which two ASEPs use the same Poisson clocks. Consider two such processes, $\eta_t$ and $\widetilde\eta_t$, started from configurations which agree everywhere except, let us say,
at the origin, where $\eta_0(0)=0$ and $\widetilde\eta_0(0)=1$. Under this coupling, the two configurations differ at exactly one site at each time $t\geq 0$, and we denote the position of this discrepancy by $Q(t)$. This discrepancy is called a \textit{second-class particle}. It behaves as a particle in interactions with holes and as a hole in interactions with ordinary (first-class) particles. More precisely, writing $0$, $1$, and $2$ for a hole, a first-class particle, and the second-class particle, respectively, the swaps $(2,0)\to(0,2)$ and $(1,2)\to(2,1)$ occur at rate $1$, while the reverse swaps occur at rate $q$. In this sense, the second-class particle tracks the propagation of a local perturbation without affecting the dynamics of the first-class particles.

Second-class particles have been extensively studied
in ASEP and related interacting particle systems. We briefly recall some of the relevant results.

The macroscopic evolution of the particle density in ASEP
is described by the conservation law
$
\partial_t\rho+\partial_x j(\rho)=0,
$
where $\rho=\rho(t,x)$ denotes the macroscopic particle density.
In this context, second-class particles play the role of
microscopic characteristics.
Depending on the initial density profile, the solution may develop shocks or rarefaction fans, and the second-class particle has different macroscopic behavior in these regimes.

In the shock regime, a second-class particle can track
the position of the macroscopic discontinuity.
Early work on microscopic shocks and travelling waves includes
De~Masi--Kipnis--Presutti--Saada~\cite{DeMasiKipnisPresuttiSaada1989},
Ferrari--Kipnis--Saada~\cite{FerrariKipnisSaada1991}, and Ferrari~\cite{Ferrari1992}. Sepp\"al\"ainen~\cite{Seppalainen2001} established laws of large numbers for second-class particles
in totally asymmetric $K$-exclusion processes, extending Ferrari's result on second-class particles
following macroscopic shocks. For Bernoulli shock data, Ferrari--Fontes~\cite{FerrariFontes1994}
established Gaussian fluctuations of order $t^{1/2}$,
arising from the randomness of the initial configuration. For certain deterministic shock initial data in the totally asymmetric
case (TASEP, $q=0$), Ferrari--Ghosal--Nejjar~\cite{FerrariGhosalNejjar2019}
determined the limiting distribution of a second-class particle on the $t^{1/3}$ fluctuation scale. Ferrari--Nejjar~\cite{FerrariNejjar2024} subsequently extended this result to joint limiting distributions at multiple times.

In a rarefaction fan, a second-class particle randomly selects one of the outgoing characteristics and follows it asymptotically, as shown for TASEP by
Ferrari--Kipnis~\cite{FerrariKipnis1995}.
Mountford--Guiol~\cite{MountfordGuiol2005} subsequently strengthened this result by proving almost sure convergence. Beyond the asymptotic motion of a single particle, Ferrari--Gon\c{c}alves--Martin~\cite{FerrariGoncalvesMartin2009} computed collision probabilities for two second-class particles in the rarefaction fan of ASEP.

In stationary ASEP, second-class particles also play an important role in the study of current fluctuations and space-time correlations. Bal\'azs--Sepp\"al\"ainen~\cite{BalazsSeppalainen2007}
established exact identities relating these quantities to the distribution of a second-class particle and later proved sharp bounds of order $t^{2/3}$ for the variance
of the current along characteristics and of order $t^{1/3}$ for the diffusivity~\cite{BalazsSeppalainen2010}.

For TASEP, directed last-passage percolation provides a geometric interpretation of second-class particles.
Ferrari--Pimentel~\cite{FerrariPimentel2005} showed that the trajectory of a second-class particle can be represented by a competition interface. The asymptotic direction and fluctuations of competition
interfaces were further studied by
Ferrari--Martin--Pimentel~\cite{FerrariMartinPimentel2006}.
More recently, Rahman--Vir\'ag~\cite{RahmanVirag2025} proved that, under KPZ scaling, the trajectory
of a second-class particle converges to a competition
interface in the directed landscape.

In the case of ASEP with $q>0$, no analogue of the competition-interface representation in last-passage percolation is known. Nevertheless, Tracy--Widom~\cite{TracyWidom2009}
obtained an exact distribution formula for a special step-type initial condition.
KPZ fluctuations at a special ASEP shock were later established by Nejjar~\cite{Nejjar2020} using a mixing argument. Adams--Bal\'azs--Jay~\cite{AdamsBalazsJay2026} recently studied second-class particles in blocking ASEP.

Of particular relevance to our work is an approach to studying second-class particles that exploits the connection between ASEP and Hecke algebras, in particular the resulting color-position symmetry in multispecies ASEP, as in Bufetov~\cite{Bufetov2020} and Borodin--Bufetov~\cite{BorodinBufetov2021}. The same symmetry was recently used by Chen~\cite{Chen2026} to derive distribution formulas for second-class particles in half-line open TASEP.

The purpose of the present paper is \textit{to derive an exact
formula for the one-time distribution of a second-class particle
in full-line ASEP with an arbitrary initial condition};
as we will see, this formula is suitable for asymptotic analysis.
The key new tool is a random \emph{oriented distance} between
binary configurations, for which we establish a duality
(Theorem~\ref{thm:intro-duality}).  Using this duality, we prove Theorem~\ref{thm:intro-exact},
which is the main result of our work. We then apply this formula to obtain the limiting distribution for arbitrary initial conditions forming a rarefaction fan
(Theorem~\ref{thm:intro-rarefaction}). We also briefly discuss ongoing work on further applications of Theorem~\ref{thm:intro-exact}.

Our proofs use the framework of random walks on Hecke algebras~\cite{Bufetov2020}, building on ideas from Borodin--Bufetov~\cite{BorodinBufetov2021}, Bufetov--Nejjar~\cite{BufetovNejjar2022}, and Borodin--Bufetov--Ferrari~\cite{BorodinBufetovFerrari2024}. In the rest of the introduction, we formulate our results in detail.
	
	\subsection{Oriented distance}\label{sec:intro-duality}

On a finite interval $I=\{1,\ldots,n\}$, let $\eta,\xi\in\{0,1\}^I$.
	Let us define $H_\eta(k)=\sum_{j=k}^n\eta(j)$, with $H_\eta(n+1)=0$,
	and consider the right-to-left discrepancy path
	\[
	g_j=H_\xi(n+1-j)-H_\eta(n+1-j),\qquad 0\le j\le n.
	\]
For any path $g=(g_0,\ldots,g_n)$ starting at zero with increments in $\{-1,0,1\}$, we define its \emph{$q$-soft maximum} $\soft(g)$ by scanning the path from left to right, starting at $Z_0=0$. At an up-step $g_{j+1}=g_j+1$, we set $Z_{j+1}=Z_j+1$ with probability
$q^{Z_j-g_j}$ and $Z_{j+1}=Z_j$ otherwise.
At flat and down-steps, we set $Z_{j+1}=Z_j$. All random choices are independent, and we denote the resulting random variable $Z_n$ by $\soft(g)$. The \emph{oriented distance} between $\eta$ and $\xi$ is then defined by
\[
\distq(\eta,\xi)=\soft(g(\eta,\xi)).
\]

Let $T_i$ denote the binary ASEP bond update at sites
$i$ and $i+1$. It swaps $10\to01$ deterministically and $01\to10$ with probability $q$.
Let $U_i$ be the particle-hole-conjugate update, which swaps $01\to10$ deterministically and
$10\to01$ with probability $q$. We write the action of these updates on configurations
from the right as $\eta T_i$ and $\eta U_i$.

The crucial property (and \emph{raison d'être}) of the
oriented distance is the following duality.
\begin{theorem}[Theorem \ref{thm:duality}]\label{thm:intro-duality}
For every $\eta,\xi\in\{0,1\}^I$ and $1\le i<n$, we have
\begin{equation}\label{eq:intro-distance-duality}
			\distq(\eta T_i,\xi)
			\overset{\mathrm d}=\distq(\eta,\xi U_i).
		\end{equation}
	\end{theorem}
The statement is elementary and, nevertheless, it seems to be new.  Although it might not be immediately clear why this duality is useful, its importance can be seen from the proofs of our next results; in short, it seems that it can be combined well with the technique of random walks on Hecke algebras. 
	\subsection{Exact finite-time formula}\label{sec:intro-exact}

	Fix an arbitrary deterministic background
	$\zeta\in\{0,1\}^{\mathbb Z\setminus\{0\}}$, and let $Q_\zeta(t)$ be
	the position of a second-class particle placed at the origin.
	For $x\in\mathbb Z$, let $\eta^x$ be an (independent) single-species ASEP
	started from the shifted step initial configuration $\eta_0^x(y)=\ind_{\{y\le x\}}$. For a binary configuration $\eta$, introduce two paths read outward from the origin, defined for $j\ge0$ by
	\begin{align}
		G_j^-(\eta)&=\sum_{y=-j}^{-1}(1-\zeta(y)-\eta(y)),
		\label{eq:intro-exact-left}\\
		G_j^+(\eta)&=\sum_{y=1}^{j}(\eta(y)+\zeta(y)-1).
		\label{eq:intro-exact-right}
	\end{align}
Both paths start at zero. At time $t$, for each $M\ge0$, we apply the $q$-soft-maximum scan to the first $M$ steps of each path, coupling the scans for different $M$ so that
each extends the preceding one. The left and right scans are independent conditionally on $\eta_t^x$. We then define

\[
D^\pm=\lim_{M\to\infty}
\soft\bigl((G_j^\pm(\eta_t^x))_{j=0}^M\bigr),
\qquad C=\eta_t^x(0).
\]
	The limits are finite almost surely: the shifted-step ASEP has only
	finitely many defects from its initial step at each finite time, and
	both paths have eventually no up-steps.

For the formula, we use the following notation
	\begin{equation}\label{eq:intro-F-kernel}
		(a;q)_m=\prod_{j=0}^{m-1}(1-aq^j),\qquad
		F_h(x)=\sum_{a=0}^h x^a(q^{-h};q)_a,
	\end{equation}
	and, for $h,r\in\mathbb Z_{\ge0}$ and $c\in\{0,1\}$, set
	\begin{equation}\label{eq:intro-kernel}
		K_q(h,r,c)=
		\begin{cases}
			F_h(q^{r+1}),&c=1,\\[2pt]
			(1-q^r)F_h(q^r),&c=0.
		\end{cases}
	\end{equation}
Our main exact result is the following.
	\begin{theorem}[Theorem \ref{thm:second-class}]\label{thm:intro-exact}
		 For every $t\ge0$ and $x\in\mathbb Z$, we have
		\begin{equation}\label{eq:intro-exact-cdf}
			\displaystyle
			\PP\bigl(Q_\zeta(t)\le x\bigr)
			=\EE\,K_q(D^-,D^+,C).
		\end{equation}
If the initial configuration $\zeta$ is random, the same identity holds after averaging over $\zeta$.
	\end{theorem}
This formula can be viewed as an analogue of the competition-interface representation of Ferrari--Pimentel~\cite{FerrariPimentel2005} for TASEP. The advantage of this representation is that the auxiliary ASEP starts from a shifted step configuration, which has been extensively studied, while the dependence on the initial background is encoded in the two discrepancy paths.

As a first asymptotic application, we use this formula to determine the limiting distribution in a rarefaction fan.
\subsection{Rarefaction fan}\label{sec:intro-rarefaction}
	
	Let $\zeta\in\{0,1\}^{\mathbb Z\setminus\{0\}}$ be a deterministic
	initial condition and place a second-class particle at the origin.  Assume that
	the one-sided empirical densities exist and satisfy
	\begin{equation}\label{eq:intro-rarefaction-densities}
		\frac1n\sum_{j=1}^n\zeta(-j)\longrightarrow p_L,
		\qquad
		\frac1n\sum_{j=1}^n\zeta(j)\longrightarrow p_R,
		\qquad 0\le p_R<p_L\le1.
	\end{equation}
	The two fan edges are
	\begin{equation}\label{eq:intro-rarefaction-edges}
		v_L=\gamma(1-2p_L),\qquad v_R=\gamma(1-2p_R),
	\end{equation}
	and for $v\in(v_L,v_R)$ we set
	\begin{equation}\label{eq:intro-rarefaction-rho}
		\rho(v)=\frac12\left(1-\frac v\gamma\right)\in(p_R,p_L).
	\end{equation}
Given $\rho\in(p_R,p_L)$, let $(B_j)_{j\in\mathbb Z}$ be i.i.d.
Bernoulli$(\rho)$ random variables, set $S_0^\pm(\rho)=0$, and define
	\begin{equation}\label{eq:intro-static-paths}
		S_n^-(\rho)=\sum_{j=1}^n(B_{-j}-\zeta(-j)),\qquad
		S_n^+(\rho)=\sum_{j=1}^n(\zeta(j)-B_j),\qquad n\ge1.
	\end{equation}
We construct the soft maxima of finite prefixes as above, using independent scans for the two paths, and define 
\begin{equation}\label{eq:intro-static-prefix-distances}
		H_{\rho,n}=\soft\bigl((S_k^-(\rho))_{k=0}^n\bigr),\qquad
		R_{\rho,n}=\soft\bigl((S_k^+(\rho))_{k=0}^n\bigr).
	\end{equation}
Since $S_n^-(\rho)/n\to\rho-p_L<0$ and
$S_n^+(\rho)/n\to p_R-\rho<0$ almost surely,
the two prefix scans stabilize almost surely
(see Lemma~\ref{lem:fan-estimates}). Thus, the a.s. limits
	\begin{equation}\label{eq:intro-static-distances}
		H_\rho=\lim_{n\to\infty}H_{\rho,n},\qquad
		R_\rho=\lim_{n\to\infty}R_{\rho,n},
	\end{equation}
	exist. The random variables $H_\rho$, $R_\rho$,
	and $B_0$ are independent.
		
Finally, we define the following function of the
initial background $\zeta$:
\begin{equation}\label{eq:intro-Phi}
    \Phi_{q,\zeta}(\rho)
    =\mathbb E\,K_q(H_\rho,R_\rho,1-B_0),
    \qquad p_R<\rho<p_L.
\end{equation}

	\begin{theorem}
    [Theorem~\ref{thm:general-rarefaction}]\label{thm:intro-rarefaction}
		Under assumption~\eqref{eq:intro-rarefaction-densities}
and with the notation above, we have
		\begin{equation}\label{eq:intro-rarefaction-limit}
			\frac{Q_\zeta(t)}{t}\Longrightarrow V_\zeta,
			\qquad t\to\infty,
		\end{equation}
where $V_\zeta$ is a random variable supported on $[v_L,v_R]$ whose distribution function is given by
		\begin{equation}\label{eq:intro-rarefaction-cdf}
        \mathbb P(V_\zeta\le v)
        =\Phi_{q,\zeta}(\rho(v)),
        \qquad v_L<v<v_R.
    \end{equation}
	\end{theorem}	

For TASEP, Ferrari--Kipnis~\cite{FerrariKipnis1995}
established the limiting speed distribution of a second-class
particle for two-sided Bernoulli rarefaction data.
Amir--Angel--Valk\'o~\cite{AmirAngelValko2011}
studied the joint asymptotic speeds of particles in
multispecies TASEP and introduced the TASEP speed process.
Cator--Pimentel~\cite{CatorPimentel2013}
extended the result of Ferrari--Kipnis to arbitrary
deterministic initial conditions forming a rarefaction fan,
using the competition-interface representation and
Busemann functions.

For ASEP with $q>0$, fewer results are available for general
initial conditions.  Ferrari--Gon\c{c}alves--Martin~
\cite{FerrariGoncalvesMartin2009} obtained the uniform limiting speed distribution for step and two-sided Bernoulli initial conditions. Borodin--Bufetov~\cite{BorodinBufetov2021} obtained
limiting distributions for finite perturbations of the
step initial condition. Aggarwal--Corwin--Ghosal~
\cite{AggarwalCorwinGhosal2023} proved the almost sure existence of the asymptotic speed for step initial data and constructed the ASEP speed process. Theorem~\ref{thm:intro-rarefaction} extends the result
of Cator--Pimentel to ASEP with $q>0$.

In Section~\ref{sec:rarefaction-examples}, we evaluate $\Phi_{q,\zeta}$ explicitly and determine the corresponding limiting distributions for various classes of initial conditions. In particular, we recover the uniform limiting distribution for two-sided Bernoulli data, derive an explicit formula for a fairly general class of initial configurations with an empty right half-line, and compute the limiting distribution for a one-site perturbation of a Bernoulli background.

\textbf{Further work.}
Theorem~\ref{thm:intro-exact} is also suitable for asymptotic analysis in several other regimes. In ongoing work~\cite{BufetovDordzhievKPZ}, we apply this theorem to obtain limiting distributions of a second-class particle in the equal-density and shock regimes for ASEP with fixed $0<q<1$, as well as in the
critical weak-asymmetry regime $q\to1$. In separate work~\cite{BufetovDordzhievTASEP}, we obtain an analogous formula for the joint distribution of a second-class and a third-class particle in TASEP ($q=0$) and apply it to several asymptotic regimes.

\textbf{AI use.}
Main results and ideas of the paper were developed by the authors without the use of AI. The writing of the paper was done with the use of ChatGPT. 

\textbf{Acknowledgements.}
A. Bufetov was partially supported by the European Research Council (ERC), Grant Agreement No. 101041499. A. Dordzhiev thanks the Max Planck Institute for Mathematics in the Sciences for providing an excellent research environment.

\section{Preliminaries}\label{sec:preliminaries}	
	\subsection{Hecke algebra}
	
	Let $S_n$ be the symmetric group, with adjacent transpositions
	$s_i=(i,i+1)$ and Coxeter length
	$\ell(w)=\#\{i<j:w(i)>w(j)\}$. We use the probabilistic
	normalization of the Hecke algebra (see, e.g. \cite{Bufetov2020}).
	
	\begin{definition}\label{def:hecke}
		The algebra $\HH(S_n)$ has basis $\{\T_w:w\in S_n\}$ and generators
		$\T_i=\T_{s_i}$ satisfying
		\begin{equation}\label{eq:hecke-relations}
			\T_i^2=(1-q)\T_i+q\mathbf 1,\qquad
			\T_i\T_j=\T_j\T_i\quad (|i-j|>1),\qquad
			\T_i\T_{i+1}\T_i=\T_{i+1}\T_i\T_{i+1}.
		\end{equation}
		For a reduced expression $w=s_{i_1}\cdots s_{i_m}$, set
		$\T_w=\T_{i_1}\cdots\T_{i_m}$. In particular,
		\begin{equation}\label{eq:right-hecke}
			\T_w\T_i=
			\begin{cases}
				\T_{ws_i},&\ell(ws_i)=\ell(w)+1,\\
				(1-q)\T_w+q\T_{ws_i},&\ell(ws_i)=\ell(w)-1.
			\end{cases}
		\end{equation}
\end{definition}
We write $\HH(S_I)$ for the same algebra on an
    interval $I$, and regard it as a subalgebra on
    any larger interval by extending permutations
    as the identity.

    An element $A=\sum_w a_w\T_w$ is called
    \emph{probabilistic} if $a_w\ge0$ and
    $\sum_w a_w=1$. Products and convex combinations
    of probabilistic elements are probabilistic.
    Such an element can be viewed as a random permutation
whose distribution is given by the coefficients $a_w$.

The following proposition is standard.
	\begin{proposition}\label{prop:involution}
		The linear map
		\begin{equation}\label{eq:involution}
			\iota(\T_w)=\T_{w^{-1}}
		\end{equation}
		is an involutive anti-automorphism:
		$\iota(AB)=\iota(B)\iota(A)$ and $\iota^2=\id$.
	\end{proposition}
	
\subsection{Colored and single-species ASEP}
On a finite interval $I$ of $\mathbb Z$, a
\emph{fully colored configuration} is a bijection
$\sigma:I\to I$, where $\sigma(i)$ is interpreted
as the color of the particle at site $i$.
The generator $\T_i$ swaps the colors at sites
$i$ and $i+1$ with probability $1$ if
$\sigma(i)<\sigma(i+1)$, and with probability
$q$ if $\sigma(i)>\sigma(i+1)$.
Otherwise, the colors remain unchanged.
This corresponds to the right Hecke action
\eqref{eq:right-hecke}.

For an integer threshold $z$, the projection
onto a single species is
\begin{equation}\label{eq:threshold}
    \mathsf p_z(\sigma)(i)
    =\ind_{\{\sigma(i)\le z\}}.
\end{equation}
Under this projection, $\T_i$ induces the
following update of binary configurations,
which we denote by $T_i$:
\begin{equation}\label{eq:binary-update}
    10\xrightarrow{\ T_i\ }01,\qquad
    01\xrightarrow{\ T_i\ }
    \begin{cases}
        10,&\text{with probability }q,\\
        01,&\text{with probability }1-q.
    \end{cases}
\end{equation}
The pairs $00$ and $11$ are unchanged. We also use the particle-hole-conjugate
updates
\begin{equation}\label{eq:dual-update}
    U_i=\comp T_i\comp,\qquad
    \comp\eta:=1-\eta.
\end{equation}
Here $U_i$ swaps $01$ to $10$ with probability
$1$ and $10$ to $01$ with probability $q$.
The other two pairs are unchanged.

We write $\eta T_i$ and $\eta U_i$ for the random
configurations obtained after the corresponding updates.
For $w\in S_I$, we similarly write $\eta T_w$ and
$\eta U_w$ for the random configurations associated
with the Hecke element $\T_w$.
More generally, for a probabilistic Hecke element
$A=\sum_w a_w\T_w$, we define $\eta T_A$ and
$\eta U_A$ as random configurations whose
distributions are the corresponding mixtures
with weights $a_w$.

The continuous-time dynamics is constructed using
independent rate-one Poisson clocks on each bond
$(i,i+1)$ with $i,i+1\in I$. Whenever a clock
rings, the corresponding colored update is applied.
For $t\geq 0$, let $\mathbf A_t^I$ denote the
probabilistic Hecke algebra element obtained by
averaging the ordered products of generators
corresponding to the clock rings up to time $t$. 
Since $\iota(\T_i)=\T_i$ and $\iota$ reverses
the order of products, while the reversed
sequences of updates have the same distribution,
we have $\iota(\mathbf A_t^I)=\mathbf A_t^I$. The classical Harris construction allows us to extend
the dynamics to $\mathbb Z$.

A single second-class particle can be obtained from
the multi-species ASEP through the following projection.
Fix a color $\kappa\in I$ and project colors smaller than
$\kappa$ to first-class particles $1$, color $\kappa$
to a second-class particle $2$, and larger colors
to holes $0$. The priority order is $1\prec 2\prec 0$. Adjacent distinct classes $a$ and $b$ swap
at rate $1$ if $a\prec b$, and at rate $q$
if $b\prec a$.

\begin{proposition}\label{prop:color-position}
    Let $A$ be a probabilistic element of the Hecke algebra,
    and let $\Sigma_A$ be the random colored configuration
    obtained by applying $A$ to the identity permutation. Then
    \begin{equation}\label{eq:color-position}
        \Sigma_A^{-1}\overset{\mathrm d}=\Sigma_{\iota(A)}.
    \end{equation}
\end{proposition}
\begin{proof}
    For each $w\in S_I$, applying $\T_w$ to the
    identity permutation produces $w$ deterministically.
    Thus, if $A=\sum_w a_w\T_w$, then
    $\PP(\Sigma_A=w)=a_w$.
    Since $\iota(A)=\sum_w a_w\T_{w^{-1}}$,
    the claim follows.
\end{proof}

This form of color--position symmetry from
\cite{BorodinBufetov2021,Bufetov2020}
will be used below.

\subsection{The longest word}

For an interval $I=[\ell,r]$, let $w_I$ be its
longest permutation, given by
\begin{equation}\label{eq:longest}
    w_I(x)=\ell+r-x,\qquad x\in I.
\end{equation}
On $I=\{1,\ldots,n\}$, one reduced expression is
\begin{equation}\label{eq:longest-word}
    \T_{w_I}
    =(\T_1\cdots\T_{n-1})
     (\T_1\cdots\T_{n-2})\cdots\T_1.
\end{equation}
Since $w_I^{-1}=w_I$, we have
$\iota(\T_{w_I})=\T_{w_I}$.
Applied to a strictly increasing colored word,
$\T_{w_I}$ reverses it deterministically.
For colored words that are not strictly increasing,
the resulting configuration may be random.

For binary configurations on $\{1,\ldots,n\}$, set
\begin{equation}\label{eq:packed}
    \Theta_b:=1^b0^{n-b},\qquad
    \Phi_b:=0^{n-b}1^b.
\end{equation}
Then the longest word gives the deterministic
sorting identities
\begin{equation}\label{eq:packed-longest}
    \Theta_b T_{w_I}=\Phi_b,\qquad
    \Phi_b U_{w_I}=\Theta_b.
\end{equation}
\section{ASEP oriented distance}\label{sec:distance}
In this section, we introduce the key new notion of
the paper --- the \textit{ASEP oriented distance} ---
and establish its main properties.
It is a random integer-valued variable defined
through the discrepancy path between two binary
configurations.
	\subsection{Discrepancy paths and the sampling rule}
Let $\eta,\xi\in\{0,1\}^{\{1,\ldots,n\}}$.
Define their suffix heights and the associated
right-to-left discrepancy path by
\begin{align}    
H_\eta(k)&=\sum_{x=k}^n\eta(x),
\qquad H_\eta(n+1)=0,
\label{eq:suffix-height}
\\
    g_j(\eta,\xi)&=H_\xi(n+1-j)-H_\eta(n+1-j),
    \qquad 0\le j\le n.
    \label{eq:discrepancy-path}
\end{align}
Thus $g_0=0$, the increments belong to $\{-1,0,1\}$,
and $g_n=|\xi|-|\eta|$, where
$|\eta|=\sum_{x=1}^n\eta(x)$.

For any such discrepancy path $g=(g_0,\ldots,g_n)$, define
\begin{equation}\label{eq:path-transform}
    \up(g)=\{j:0\le j<n,\ g_{j+1}=g_j+1\}.
\end{equation}
For $g=g(\eta,\xi)$, these are precisely the steps
corresponding to sites where $\eta=0$ and $\xi=1$.

\begin{definition}\label{def:distance}
     Let $g=(g_0,\ldots,g_n)$ be a discrepancy path.
    Start at $Z_0=0$ and scan $g$ from $j=0$ to $j=n-1$. At an up-step $g_j\to g_j+1$,
    set
    \begin{equation}\label{eq:soft-scan}
        Z_{j+1}=
        \begin{cases}
            Z_j+1,&\text{with probability }q^{Z_j-g_j},\\
            Z_j,&\text{with probability }1-q^{Z_j-g_j}.
        \end{cases}
    \end{equation}
    If the path stays flat or steps down, set
    $Z_{j+1}=Z_j$. At each up-step, the choices are independent.
    We call the resulting random variable
    $\soft(g)=Z_n$ the \emph{$q$-soft maximum} of $g$. For two binary configurations $\eta,\xi$,
    the \emph{ASEP oriented distance} is the
    random variable
\begin{equation}\label{eq:distance}
\distq(\eta,\xi)=\soft\bigl(g(\eta,\xi)\bigr).
    \end{equation}
\end{definition}
To ensure that the construction is well-defined, we need to verify that \(q^{Z_j-g_j}\leq1\) at every up-step. 
	\begin{lemma}\label{lem:scan-validity}    
For every $0\le j\le n$, we have
		\begin{equation}\label{eq:running-record}
			Z_j\ge\max_{0\le i\le j}g_i.
		\end{equation}
		In particular, the procedure above is well-defined and
		$0\le\soft(g)\le|\up(g)|$.
	\end{lemma}
\begin{proof}
We argue by induction on $j$. The claim holds at $j=0$.
Suppose it holds at $j$.
Since $Z_{j+1}\ge Z_j$, the inequality can fail
at $j+1$ only if $g_{j+1}>Z_j$.
As $g_{j+1}\le g_j+1\le Z_j+1$,
this implies $g_j=Z_j$ and $g_{j+1}=g_j+1$.
At this up-step, $Z_{j+1}=Z_j+1$ with probability
$q^{Z_j-g_j}=1$, proving \eqref{eq:running-record}. The bounds on $\soft(g)$ follow
since $Z_j$ is nonnegative and increases only at up-steps.
\end{proof}	
\begin{remark}
The $q$-soft maximum can also be sampled by reading the path in reverse order; see
Lemma~\ref{lem:reverse-queue} in Section~\ref{sec:rarefaction-examples}. We note that the transition probabilities of this reverse-order algorithm also appear in Martin's multi-line queue construction~\cite{Martin2020}.

\end{remark}

\subsection{The polynomial transform}
	
	For $z\in\NN$ and a discrepancy path $g$, define
\begin{equation}\label{eq:packed-polynomial}
    \begin{aligned}
        \poly_z(u)&=\prod_{a=0}^{z-1}(1+u q^{-a}),
        \qquad \poly_0(u)=1,\\
        \Pi_g(u)&=\prod_{j\in\up(g)}(1+u q^{-g_j}).
    \end{aligned}
\end{equation}
For binary configurations $\eta,\xi$, write
\begin{equation}\label{eq:Phi}
\Pi_u(\eta,\xi):=\Pi_{g(\eta,\xi)}(u).
\end{equation}
    \begin{proposition}\label{prop:transform}
		Let \(g\) be a finite discrepancy path. Then
		\begin{equation}\label{eq:transform}
			\EE\,\poly_{\soft(g)}(u)=\Pi_g(u).
		\end{equation}
		This identity uniquely determines the distribution
    of $\soft(g)$. In particular, one has $\EE\,\poly_{\distq(\eta,\xi)}(u)=\Pi_u(\eta,\xi)$.
	\end{proposition}
	
	\begin{proof} 
If the current state is \(Z_j=z\) and the next up-step starts
at height \(y\), then
\[
\begin{aligned}
\EE[\poly_{Z_{j+1}}(u)\mid Z_j=z]
&=(1-q^{z-y})\poly_z(u)+q^{z-y}\poly_{z+1}(u)\\
&=\poly_z(u)(1+uq^{-y}).
\end{aligned}
\]
Iteration gives \eqref{eq:transform}. Uniqueness follows since the polynomials
$\poly_0,\ldots,\poly_m$ have distinct degrees and nonzero leading coefficients, and hence form a basis.

\end{proof}

	\begin{example}\label{ex:excursion}
		For $\eta=01$ and $\xi=10$, the discrepancy path \(g(\eta,\xi)\) is
		$(0,-1,0)$, and $\distq(\eta,\xi)$ is Bernoulli with parameter $q$. Interchanging $\eta$ and $\xi$, the discrepancy
path becomes $(0,1,0)$, and
$\distq(\xi,\eta)=1$ deterministically.
 Note that our oriented distance is not symmetric.
	\end{example}

\begin{example}\label{ex:distance-two-rises}
    For $\eta=001$ and $\xi=110$, the discrepancy
    path is $g=(0,-1,0,1)$, with two up-steps
    starting at heights $-1$ and $0$.
    At the first up-step, the scan increases with
    probability $q$. At the second, the increase
    is forced if the first one failed, and occurs
    with probability $q$ otherwise. Therefore,
    \begin{equation}\label{eq:two-rise-law}
        \PP(\distq(001,110)=1)=1-q^2,\qquad
        \PP(\distq(001,110)=2)=q^2.
    \end{equation}
\end{example}

\begin{example}\label{ex:packed-distances}
    Let $\eta\in\{0,1\}^n$ be arbitrary, and let
    $
        \Theta_b=1^b0^{n-b},
        \Phi_b=0^{n-b}1^b.
    $
    We have
    \[
        \distq(\eta,\Phi_b)
        =\sum_{i=n-b+1}^n(1-\eta(i)),
        \qquad
        \distq(\Theta_b,\eta)
        =\sum_{i=b+1}^n\eta(i).
    \]
 Thus, the first distance equals the
    number of holes in the rightmost $b$ sites of $\eta$,
    while the second equals the number of particles
    in the rightmost $n-b$ sites. Indeed, in both cases, all up-steps of the discrepancy
    path precede all down-steps, so the soft maximum
    increases at every up-step.

    In contrast, $\distq(\eta,\Theta_b)$ and
    $\distq(\Phi_b,\eta)$ are generally random, since
    their discrepancy paths have all down-steps preceding
    all up-steps. We study their distributions in
    Section~\ref{sec:pure-death}.
\end{example}
	\subsection{Local duality}
We now establish the duality of the oriented distance,
together with its extension to probabilistic Hecke elements.

\begin{theorem}\label{thm:duality}
    For every bond $i$ and every pair of binary
    configurations $\eta,\xi$, we have
    \begin{equation}\label{eq:distance-duality}
        \distq(\eta T_i,\xi)\overset{\mathrm d}=
        \distq(\eta,\xi U_i).
    \end{equation}
    More generally, for every probabilistic Hecke
    element $A$,
    \begin{equation}\label{eq:word-duality}
        \distq(\eta T_A,\xi)\overset{\mathrm d}=
        \distq(\eta,\xi U_{\iota(A)}).
    \end{equation}
\end{theorem}

\begin{proof}
    
By Proposition~\ref{prop:transform}, it suffices to establish the polynomial identity $\EE\,\Pi_u(\eta T_i,\xi)
=\EE\,\Pi_u(\eta,\xi U_i)$. Write $g=g(\eta,\xi)$ for the discrepancy path. Only the middle height $g_{n-i}$ of the triple
$(g_{n-i-1},g_{n-i},g_{n-i+1})$ can change
under either update.
Since $\Pi_g(u)$ is a product over the up-steps
of $g$, it suffices to compare the local factors.

    An $\eta$-pair $10$ or a $\xi$-pair $01$
    lowers this height by one deterministically,
    while an $\eta$-pair $01$ or a $\xi$-pair $10$
    raises it with probability $q$.

    If both pairs are constant, neither update changes
    the path. If exactly one pair is nonconstant,
    its update either moves an up-step to the neighboring
    edge without changing its starting height or
    leaves the local factors unchanged.
    If the two nonconstant pairs are oppositely ordered,
    the updates induce the same random modification
    of the path.
    
    The remaining cases are $10,10$ and $01,01$.
    In both cases, the local height triple is $(y,y,y)$.
    Lowering its middle height creates an up-step
    with factor $1+uq^{1-y}$, whereas raising it
    creates an up-step with factor $1+uq^{-y}$.
    The identity
    \[1+uq^{1-y}=(1-q)+q(1+uq^{-y})\]
    gives the required equality of expected local factors.
    This proves the polynomial identity.

    By Proposition~\ref{prop:transform}, the polynomial
    identity implies \eqref{eq:distance-duality}.
    Iterating the local duality transfers successive
    updates to the second argument in reverse order.
    Taking mixtures gives \eqref{eq:word-duality}.
\end{proof}

\begin{example}
    Let $\eta=\xi=10$. The update $T_1$ sends
    $\eta$ to $01$ deterministically. Hence,
    $\distq(\eta T_1,\xi)$ is Bernoulli with
    parameter $q$ by Example~\ref{ex:excursion}. On the other hand, $\xi U_1$ equals $01$
    with probability $q$ and $10$ with probability
    $1-q$. Since $\distq(\eta,01)=1$ and
$\distq(\eta,10)=0$ deterministically,
$\distq(\eta,\xi U_1)$ is also Bernoulli
with parameter $q$.

\end{example}

	\section{The origin after the longest word}\label{sec:longest}
In this section, we study the action of the longest-word element $T_{w_I}$ on binary configurations on a finite interval $I$ and compute the resulting one-point occupation probabilities.	
\subsection{Packed-window transforms}
    Let $I=\{1,\ldots,n\}$, and let
$X=\eta T_{w_I}$, where \(\eta\) is an arbitrary binary configuration on $I$. For $0\le b\le n$, let $H_b(\omega)$ denote the
number of holes in the rightmost $b$ sites of a binary
configuration $\omega$, with $H_0(\omega)=0$.

\begin{lemma}\label{lem:packed-window}
    For every $0\le b\le n$,
    \begin{equation}\label{eq:packed-window}
        H_b(X)\overset{\mathrm d}=
        \distq(\eta,\Theta_b).
    \end{equation}
\end{lemma}
\begin{proof}
    By Example~\ref{ex:packed-distances}, we have
    $H_b(X)=\distq(X,\Phi_b)$.
Applying Theorem~\ref{thm:duality},
    together with $\iota(\T_{w_I})=\T_{w_I}$
    and $\Phi_b U_{w_I}=\Theta_b$, gives
    \[
        \distq(X,\Phi_b)
        \overset{\mathrm d}=
        \distq(\eta,\Phi_b U_{w_I})
        =\distq(\eta,\Theta_b).
    \]
\end{proof}

\subsection{A pure-death chain}\label{sec:pure-death}

We will need the following auxiliary object. Let $(Y_s)_{s\ge0}$ be a discrete-time Markov chain on
$\mathbb Z_{\ge0}$ with transitions
\begin{equation}\label{eq:death-chain}
    Y\longmapsto
    \begin{cases}
        Y,&\text{with probability }q^Y,\\
        Y-1,&\text{with probability }1-q^Y.
    \end{cases}
\end{equation}

\begin{lemma}\label{lem:packed-distance-death}
    Let $\eta\in\{0,1\}^n$ and
    $\Theta_b=1^b0^{n-b}$.
    Set
    \[
        h=\sum_{i=1}^b(1-\eta(i)),
        \qquad
        r=\sum_{i=b+1}^n\eta(i),
    \]
    where $h$ is the number of holes in the leftmost
    $b$ sites of $\eta$, and $r$ is the number of
particles in the rightmost $n-b$ sites. Then
    \[
        \distq(\eta,\Theta_b)
        \overset{\mathrm d}=h-r+Y_h,
        \qquad Y_0=r.
    \]
The same argument applies to $\distq(\Phi_b,\eta)$,
with $h$ replaced by the number of particles in the leftmost
$n-b$ sites and $r$ by the number of holes in the rightmost
$b$ sites.
\end{lemma}
\begin{proof}
    Ignoring horizontal steps, the discrepancy path
    for $(\eta,\Theta_b)$ consists of $r$ down-steps
    followed by $h$ up-steps. After all $r$ down-steps, the path is at height $-r$,
while the soft maximum remains zero.
   At each of the following $h$ up-steps, the difference
between the soft maximum and the path stays unchanged
with probability $q^Y$ and decreases by one otherwise,
so it evolves according to the Markov chain
\eqref{eq:death-chain}.
Since the final height of the path is $h-r$,
we obtain $\distq(\eta,\Theta_b)\overset{\mathrm d}=h-r+Y_h$.
\end{proof}
For $r,s\ge0$, define
\[
    m_s(r)=\EE[Y_s\mid Y_0=r],
    \qquad
    a_s(r)=m_s(r+1)-m_s(r).
\]

\begin{lemma}\label{lem:death-chain}
    We have $a_0(r)=1$ and the recursion
    \begin{equation}\label{eq:positive-kernel-recursion}
        \begin{aligned}
            a_{s+1}(0)&=q\,a_s(0),\\
            a_{s+1}(r)&=q^{r+1}a_s(r)
                +(1-q^r)a_s(r-1),
                &&r\ge1.
        \end{aligned}
    \end{equation}
    The recursion admits the explicit solution
    \begin{equation}\label{eq:F}
    a_s(r)=F_s(q^{r+1}),\qquad
    F_s(x)=\sum_{j=0}^{s}x^j(q^{-s};q)_j,
\end{equation}
where $(v;q)_j=\prod_{i=0}^{j-1}(1-vq^i)$
is the $q$-Pochhammer symbol. In particular, $a_s(0)=F_s(q)=q^s$.
\end{lemma}

\begin{proof}
    Conditioning on the first step gives
    \[
        m_{s+1}(r)=q^r m_s(r)
            +(1-q^r)m_s(r-1),\qquad r\ge1.
    \]
    Taking adjacent differences and using $m_s(0)=0$
    yields \eqref{eq:positive-kernel-recursion}.

To verify the formula, note that
the polynomials $F_s$ satisfy
\begin{equation}\label{eq:F-positive-recursion}
    F_{s+1}(qx)=qxF_s(qx)+(1-x)F_s(x).
\end{equation}
Indeed, for $1\le j\le s+1$, the coefficient of $x^j$
on the right-hand side is
\[
    (q^{-s};q)_j+(q^j-1)(q^{-s};q)_{j-1}
    =q^j(q^{-s-1};q)_j,
\]
as required.

Substituting $x=q^r$ for $r\ge0$ gives
\eqref{eq:positive-kernel-recursion}. Since $F_0=1$, the formula follows by induction,
and \eqref{eq:positive-kernel-recursion} gives
$F_s(q)=a_s(0)=q^s$.

\end{proof}
\subsection{The occupation kernel} 
\begin{definition}\label{def:kernel}
    For $h,r\in\NN$ and $c\in\{0,1\}$, set
    \begin{equation}\label{eq:kernel}
        K_q(h,r,c)=
        \begin{cases}
            F_h(q^{r+1}),&c=1,\\[2pt]
            (1-q^r)F_h(q^r),&c=0.
        \end{cases}
    \end{equation}
\end{definition}
\begin{theorem}\label{thm:longest-kernel}
    Let $I=[a,b]$, choose $k\in I$, and set
    $k^*=a+b-k$. For a binary configuration
    $\eta$ on $I$, define
    \begin{equation}\label{eq:hrc}
        h=\sum_{x=a}^{k-1}(1-\eta(x)),\qquad
        r=\sum_{x=k+1}^{b}\eta(x),\qquad
        c=\eta(k).
    \end{equation}
    Then
    \begin{equation}\label{eq:longest-kernel}
        \PP\bigl((\eta T_{w_I})(k^*)=1\bigr)
        =K_q(h,r,c).
    \end{equation}
\end{theorem}

\begin{proof}
    By relabeling, we may assume $I=\{1,\ldots,n\}$,
    so that $k^*=n+1-k$.
    Put $X=\eta T_{w_I}$. Since
    \[
        H_k(X)-H_{k-1}(X)=1-X(k^*),
    \]
    Lemma~\ref{lem:packed-window} gives
    \[
        \PP(X(k^*)=1)
        =1-\EE\distq(\eta,\Theta_k)
         +\EE\distq(\eta,\Theta_{k-1}).
    \]
For $\Theta_k$ and $\Theta_{k-1}$, the numbers
of down-steps and up-steps are $(r,h+1-c)$
and $(r+c,h)$, respectively.
Thus, Lemma~\ref{lem:packed-distance-death} gives
\[
    \EE\distq(\eta,\Theta_k)=h+1-c-r+m_{h+1-c}(r),
    \qquad
    \EE\distq(\eta,\Theta_{k-1})=h-r-c+m_h(r+c).
\]
Consequently,
\[
    \PP(X(k^*)=1)
    =m_h(r+c)-m_{h+1-c}(r).
\]
If $c=1$, Lemma~\ref{lem:death-chain} gives
    \[
        m_h(r+1)-m_h(r)
        =a_h(r)=F_h(q^{r+1}).
    \]
    If $c=0$ and $r\ge1$, the first-step recursion
    gives
    \[
    \begin{aligned}
        m_h(r)-m_{h+1}(r)
        &=(1-q^r)a_h(r-1)\\
        &=(1-q^r)F_h(q^r).
    \end{aligned}
    \]
    For $r=0$, both sides vanish.
    This proves \eqref{eq:longest-kernel}.
\end{proof}
We also note that Theorem~\ref{thm:longest-kernel},
together with particle-hole reflection symmetry, gives
\begin{equation}\label{eq:kernel-symmetry}
    K_q(h,r,c)+K_q(r,h,1-c)=1.
\end{equation}

	\section{The exact distribution of a second-class particle}
	\label{sec:second-class}
We now combine the oriented-distance duality with
the longest-word kernel from Section~\ref{sec:longest}
to derive our main exact result, Theorem~\ref{thm:second-class},
which gives the distribution of the position of a
second-class particle at a fixed time.

\subsection{The full-line distribution formula}

Fix a deterministic background
$\zeta\in\{0,1\}^{\ZZ\setminus\{0\}}$.
Consider ASEP with a second-class particle initially at
the origin and first-class particles and holes elsewhere
according to $\zeta$.
Let $Q_\zeta(t)$ denote the position of the second-class
particle at time $t$.

For a binary configuration $\eta$ on $\ZZ$ and $j\ge0$,
define
\begin{align}
    G_j^-(\eta)
    &=\sum_{y=-j}^{-1}
    \bigl((1-\zeta(y))-\eta(y)\bigr),
    \label{eq:left-path}\\
    G_j^+(\eta)
    &=\sum_{y=1}^{j}
    \bigl(\eta(y)-(1-\zeta(y))\bigr).
    \label{eq:right-path}
\end{align}
Both paths start at zero, with empty sums for $j=0$.
The left path is read in the order $-1,-2,\ldots$,
whereas the right path is read in the order $1,2,\ldots$.

For $M\ge0$, write
$G^{\pm,[M]}(\eta)=(G_j^\pm(\eta))_{j=0}^M$
and $D^{\pm,[M]}(\eta)=\soft(G^{\pm,[M]}(\eta))$. We couple the soft maxima for different $M$ so that each scan extends the preceding one.  The two scans are independent conditionally on $\eta$.

Fix $t\ge0$ and $x\in\ZZ$, and let $\eta^x$ be an
auxiliary ASEP started from the shifted step configuration
\begin{equation}\label{eq:shifted-step}
    \eta_0^x(y)=\ind_{\{y\le x\}}.
\end{equation}
At time $t$, the configuration $\eta_t^x$ is almost surely
equal to one sufficiently far to the left and to zero
sufficiently far to the right.
Consequently, both paths $G^-(\eta_t^x)$ and
$G^+(\eta_t^x)$ have only finitely many up-steps.
Since the soft maximum can increase only at an up-step,
the soft maxima of their finite prefixes stabilize
almost surely. We may therefore define
\begin{equation}\label{eq:full-line-distances}
D^\pm=\lim_{M\to\infty}D^{\pm,[M]}(\eta_t^x),
    \qquad C=\eta_t^x(0).
\end{equation}
The expectation below is taken over the auxiliary ASEP and the two soft-maximum scans.

\begin{theorem}\label{thm:second-class}
    For every $t\ge0$ and $x\in\ZZ$,
    \begin{equation}\label{eq:second-class-formula}
        \PP\bigl(Q_\zeta(t)\le x\bigr)
        =\EE\,K_q(D^-,D^+,C).
    \end{equation}
\end{theorem}
In particular, the dependence on the initial condition
is contained in the two discrepancy paths.
If the initial condition is random, the same identity holds after averaging over it.

The remainder of this section is devoted to the proof
of Theorem~\ref{thm:second-class}.
\subsection{Shaping packed configurations}
	
	The next lemma shows the main idea how to create an arbitrary initial condition of ASEP with the help of Hecke algebra.

\begin{lemma}\label{lem:shaping}
    Let $\xi\in\{0,1\}^{\{1,\ldots,n\}}$ have $k$
    particles at positions $p_1<\cdots<p_k$.
    There exists a word $\mathbf D_\xi$ of length
    $\sum_{j=1}^{k}(p_j-j)$ such that
    \begin{equation}\label{eq:shaping}
        \Theta_k T_{\mathbf D_\xi}=\xi,\qquad
        \Phi_{n-k}U_{\mathbf D_\xi}=\comp\xi
    \end{equation}
    deterministically. Moreover, for any binary
    configuration $\eta$ on $\{1,\ldots,n\}$,
    \begin{equation}\label{eq:shaping-defect}
        \sum_{y=k+1}^{n}
        \bigl(1-(\eta T_{\iota(\mathbf D_\xi)})(y)\bigr)
        \overset{\mathrm d}=\distq(\eta,\comp\xi).
    \end{equation}
    The identity \eqref{eq:shaping-defect} holds for
    any word that sends $\Theta_k$ to $\xi$
    deterministically.
\end{lemma}
	
	\begin{proof}
		Starting from $\Theta_k=1^k0^{n-k}$, place the particles in the
		order $k,k-1,\ldots,1$. To place particle $j$, apply the generators
$T_j,T_{j+1},\ldots,T_{p_j-1}$ in this order. Every swap has the form $10\mapsto01$ and therefore
occurs deterministically.
The resulting word $\mathbf D_\xi$ has length $\sum_{j=1}^k(p_j-j)$. The second identity in \eqref{eq:shaping} follows by particle-hole complementation, since $\comp\Theta_k=\Phi_{n-k}$ and $U_i=\comp T_i\comp$.

By Example~\ref{ex:packed-distances}, the oriented
distance to $\Phi_{n-k}$ is deterministic and equals to
the number of holes in the rightmost $n-k$ sites.
Hence,
\[
    \sum_{y=k+1}^{n}
    \bigl(1-(\eta T_{\iota(\mathbf D_\xi)})(y)\bigr)
    =\distq\bigl(\eta T_{\iota(\mathbf D_\xi)},
    \Phi_{n-k}\bigr).
\]
Applying the duality \eqref{eq:word-duality}, we obtain
\[
    \distq\bigl(\eta T_{\iota(\mathbf D_\xi)},
    \Phi_{n-k}\bigr)
    \overset{\mathrm d}=
    \distq\bigl(\eta,\Phi_{n-k}U_{\mathbf D_\xi}\bigr)
    =\distq(\eta,\comp\xi).
\]
	\end{proof}
	
	\subsection{Proof of the distribution formula}
	\begin{proof}[Proof of Theorem~\ref{thm:second-class}]
We first construct $\zeta$ as a projection of a colored
configuration and apply color-position symmetry to
study $Q_\zeta(t)$.

First assume that $\zeta$ is a finite perturbation
of the step configuration. Choose $M>|x|$ sufficiently
large so that
\begin{equation}\label{eq:background-tails}
    \zeta(y)=1\quad(y<-M),\qquad
    \zeta(y)=0\quad(y>M).
\end{equation}
Write $I^-=[-M,-1]$, $I^+=[1,M]$, and
$\zeta^\pm=\zeta|_{I^\pm}$. Let $a$ be the number
of holes in $I^-$ and $b$ the number of particles
in $I^+$
\begin{equation}\label{eq:charge}
    a=\sum_{y=-M}^{-1}(1-\zeta(y)),\qquad
    b=\sum_{y=1}^{M}\zeta(y),\qquad
    \kappa=b-a.
\end{equation}
By Lemma~\ref{lem:shaping}, there exist words
$\mathbf D^-$ and $\mathbf D^+$, supported on
$I^-$ and $I^+$, respectively, which send
$1^{M-a}0^a$ to $\zeta^-$ and $1^b0^{M-b}$
to $\zeta^+$.
The two words commute, and neither acts on the origin. Let us also define
\begin{equation}\label{eq:asymmetric-interval}
    I_0=[-a,b],\qquad w=w_{I_0}.
\end{equation}
Note that $w(0)=\kappa$ and $w(\kappa)=0$.

We are now in a position to construct the prescribed
initial configuration. Consider the identity coloring and apply $\T_w$.
On $I_0$, the resulting color at site $y$ is
$w(y)=\kappa-y$, while outside $I_0$ the colors
remain unchanged. Thus, color $\kappa$ is at the
origin, with larger colors to its left and smaller
colors to its right within $I_0$. Projecting colors smaller than $\kappa$ to first-class
particles, color $\kappa$ to a second-class particle,
and larger colors to holes gives
\[
    \bigl(1^{M-a}0^a,\mathsf 2,1^b0^{M-b}\bigr)
\]
on $I^-\cup\{0\}\cup I^+$. Applying the shaping words $\mathbf D^-$ and
$\mathbf D^+$ then produces the prescribed initial
configuration with background $\zeta$. Set $\Lambda_M=[-M,M]$ and consider the probabilistic
Hecke elements
\[
    B=\T_w\mathbf D^+\mathbf D^-\mathbf A_t^{\Lambda_M},
    \qquad
    A=\iota(B)=\mathbf A_t^{\Lambda_M}
    \iota(\mathbf D^-)\iota(\mathbf D^+)\T_w,
\]
where we used $\iota(\mathbf A_t^{\Lambda_M})
=\mathbf A_t^{\Lambda_M}$.
Let $\Sigma_B$ and $\Sigma_A$ denote the random colored
configurations with distributions given by the
coefficients of $B$ and $A$, respectively. 
By construction, $\Sigma_B^{-1}(\kappa)$ has the same
distribution as the position $Q_\zeta^M(t)$ of the
second-class particle in ASEP on $\Lambda_M$. Applying color-position symmetry
(Proposition~\ref{prop:color-position}), we obtain
\[
    \PP\bigl(Q_\zeta^M(t)\le x\bigr)
    =\PP\bigl(\Sigma_B^{-1}(\kappa)\le x\bigr)
    =\PP\bigl(\Sigma_A(\kappa)\le x\bigr).
\]
Under the projection that identifies colors
$\le x$ with particles and colors $>x$ with holes,
the identity coloring becomes the shifted step
configuration. Let $\eta^{x,M}$ denote ASEP on
$\Lambda_M$ started from this configuration. Hence,
\[
    \PP\bigl(Q_\zeta^M(t)\le x\bigr)
    =\PP\left(
    \bigl(\eta_t^{x,M}T_{\iota(\mathbf D^-)}
    T_{\iota(\mathbf D^+)}T_w\bigr)(\kappa)=1
    \right).
\]		

We now study the occupation probability on the
right-hand side using the oriented-distance duality
(Theorem~\ref{thm:duality}).
Condition on $\eta=\eta_t^{x,M}$ and write
$\widetilde\eta=\eta T_{\iota(\mathbf D^-)}
T_{\iota(\mathbf D^+)}$.
Define
\begin{equation}\label{eq:shaped-defects}
    H=\sum_{y=-a}^{-1}(1-\widetilde\eta(y)),\qquad
    R=\sum_{y=1}^{b}\widetilde\eta(y),\qquad
    C=\widetilde\eta(0)=\eta(0).
\end{equation}
Since the reflection of $0$ in $I_0$ is $\kappa$,
Theorem~\ref{thm:longest-kernel} gives
\[\PP((\widetilde\eta T_w)(\kappa)=1\mid\widetilde\eta)
=K_q(H,R,C).\]

By Lemma~\ref{lem:shaping}, $H$ has the same distribution
as $\distq(\eta|_{I^-},\comp\zeta^-)$, whose discrepancy
path is precisely $G^{-,[M]}(\eta)$.
 Hence, $H$ has the same distribution
as $\soft(G^{-,[M]}(\eta))$.
				
For the right interval, consider the particle-hole
reflection $(J\omega)(y)=1-\omega(M+1-y)$ on $I^+$.
It satisfies $JT_iJ=T_{M-i}$. Let
$\widehat{\mathbf D}^+$ be obtained from $\mathbf D^+$
by replacing each generator $\T_i$ with $\T_{M-i}$.
Since $J\Theta_b=\Theta_{M-b}$, this word sends
$\Theta_{M-b}$ to $J\zeta^+$ deterministically. Writing $\eta^+=\eta|_{I^+}$, we have
\[
    J\bigl(\eta^+T_{\iota(\mathbf D^+)}\bigr)
    \overset{\mathrm d}=
    (J\eta^+)T_{\iota(\widehat{\mathbf D}^+)}.
\]
Since $J$ transforms particles in the leftmost $b$
sites into holes in the rightmost $b$ sites,
Lemma~\ref{lem:shaping} implies that $R$ has the
same distribution as
$\distq(J\eta^+,\comp J\zeta^+)$.
The corresponding discrepancy path satisfies
\[
    \sum_{y=M+1-j}^{M}
    \bigl((\comp J\zeta^+)(y)-(J\eta^+)(y)\bigr)
    =G_j^+(\eta),\qquad 0\le j\le M.
\]
Hence, $R$ has the same distribution as
$\soft(G^{+,[M]}(\eta))$.

Since the two reversed shaping words act on disjoint
intervals, $H$ and $R$ are independent conditionally
on $\eta$. Taking expectations, we obtain
\begin{equation}\label{eq:second-class-box}
    \PP\bigl(Q_\zeta^M(t)\le x\bigr)
    =\EE\,K_q\left(
    \soft(G^{-,[M]}(\eta_t^{x,M})),
    \soft(G^{+,[M]}(\eta_t^{x,M})),
    \eta_t^{x,M}(0)\right).
\end{equation}
Now, we finally pass to the full line by a standard
localization argument.

For arbitrary $\zeta$, let $\zeta^{[M]}$ agree
with $\zeta$ on $\Lambda_M\setminus\{0\}$ and
have occupied left and empty right tails.
Apply \eqref{eq:second-class-box} to $\zeta^{[M]}$.

Couple the finite-interval and full-line ASEPs using the basic coupling. For fixed $t$, the positions of their second-class particles agree with probability tending to one as $M\to\infty$. Indeed, set $r=\lfloor M/2\rfloor$. If the two second-class particles have different positions at time $t$, then either at least one of them has left $[-r,r]$ before time $t$, or an influence from outside $\Lambda_M$ has reached $[-r,r]$. The first event has probability tending to zero, since each second-class particle jumps at total rate at most $2$ and must make at least $r$ jumps to leave $[-r,r]$.  The second event requires at least $M-r$ chronologically ordered clock rings, and hence a union bound over possible paths and standard Poisson estimates show that its probability also tends to zero.

It remains to pass to the limit on the right-hand side of \eqref{eq:second-class-box}. Couple $\eta^{x,M}$ and $\eta^x$ using the basic coupling, and note that, with probability tending to one, $\eta_s^x$ is identically one on $(-\infty,-M]$ and identically zero on $[M,\infty)$ throughout $[0,t]$.  On this event, the two processes agree on $\Lambda_M$, and neither discrepancy path has any up-steps
after its first $M$ steps. Hence, $D^{\pm,[M]}(\eta_t^{x,M})=D^\pm$. Consequently, the kernel in \eqref{eq:second-class-box}
converges in probability to $K_q(D^-,D^+,C)$, and since $0\le K_q\le1$, the expectations also converge, which completes the proof.

\end{proof}
	
\section{The rarefaction limit for general initial configurations}\label{sec:rarefaction}
We now use Theorem~\ref{thm:second-class} to study the large-time behavior of the second-class particle in the rarefaction regime.
\subsection{Statement of the result}

Let $\zeta\in\{0,1\}^{\ZZ\setminus\{0\}}$ be a deterministic configuration
with one-sided densities $0\le p_R<p_L\le1$ such that
\begin{equation}\label{eq:fan-density-assumption}
 \frac1n\sum_{j=1}^n\zeta(-j)\longrightarrow p_L,\qquad
 \frac1n\sum_{j=1}^n\zeta(j)\longrightarrow p_R,
 \qquad n\to\infty.
\end{equation}
As before, the second-class particle starts at the origin.
Set
\begin{equation}\label{eq:fan-speeds}
 \gamma=1-q,\qquad v_L=\gamma(1-2p_L),\qquad
 v_R=\gamma(1-2p_R),
\end{equation}
and define the rarefaction density profile
\begin{equation}\label{eq:fan-rho}
 \rho(v)=
 \begin{cases}
 1,&v\le-\gamma,\\
 \frac12(1-v/\gamma) ,&-\gamma<v<\gamma,\\
 0,&v\ge\gamma.
 \end{cases}
\end{equation}
For $\rho\in[0,1]$, let $(B_j)_{j\in\ZZ}$ be i.i.d.
Bernoulli$(\rho)$ random variables and set, for $n\ge0$,
\begin{equation}\label{eq:fan-static-paths}
 S_n^-(\rho)=\sum_{j=1}^n\bigl(B_{-j}-\zeta(-j)\bigr),\qquad
 S_n^+(\rho)=\sum_{j=1}^n\bigl(\zeta(j)-B_j\bigr).
\end{equation}
For $p_R<\rho<p_L$, both infinite paths
$S^\pm(\rho)=(S_n^\pm(\rho))_{n\ge0}$
tend to $-\infty$ linearly almost surely.
As in the previous section, we couple the soft maxima of finite prefixes so that each scan extends the preceding one. By Lemma~\ref{lem:fan-estimates}(a), they stabilize almost surely, and we denote their limiting values by
\begin{equation}\label{eq:fan-static-distances}
 H_\rho=\soft(S^-(\rho)),\qquad
 R_\rho=\soft(S^+(\rho)).
\end{equation}
The random variables $H_\rho,R_\rho,B_0$ are independent.
Define
\begin{equation}\label{eq:fan-functional}
 \Phi_{q,\zeta}(\rho)=\EE\,K_q\bigl(H_\rho,R_\rho,1-B_0\bigr).
\end{equation}
\begin{theorem}\label{thm:general-rarefaction}
 Under \eqref{eq:fan-density-assumption},
 \begin{equation}\label{eq:fan-weak-limit}
  \frac{Q_\zeta(t)}t\ \Longrightarrow\ V_\zeta,
  \qquad t\longrightarrow\infty,
 \end{equation}
where the random variable $V_\zeta$ is supported on $[v_L,v_R]$.
Its distribution function is continuous on $(v_L,v_R)$
and is given by
\begin{equation}\label{eq:fan-cdf}
 \PP(V_\zeta\le v)=\Phi_{q,\zeta}(\rho(v)),
 \qquad v_L<v<v_R.
\end{equation}
\end{theorem}
\begin{remark}
If the initial configuration $\zeta$ is random, the theorem applies
conditionally on $\zeta$. The unconditional limiting distribution
function in the fan is then
$\EE_\zeta\Phi_{q,\zeta}(\rho(v))$.
\end{remark}
The rest of this section is devoted to proving the theorem.
\subsection{Technical estimates}\label{sec:estimates}
We begin by proving several technical estimates needed for the asymptotic analysis.

\begin{lemma}\label{lem:fan-estimates}
(a) Let $g=(g_j)_{j\ge0}$ be a path starting at zero
with increments in $\{-1,0,1\}$.
Suppose that $g_j\le-\delta j$ for all $j\ge M$,
for some $\delta>0$ and $M\ge0$.
Then $\soft(g^{[N]})$ stabilizes almost surely as $N\to\infty$.
Moreover,
\begin{equation}\label{eq:fan-scan-tail}
\PP\bigl(\soft(g)\ne\soft(g^{[M]})\bigr)
 \le \frac{q^{\delta M}}{1-q^\delta}.
\end{equation}

(b) For $h,r\in\NN$ and $c\in\{0,1\}$, set
$d=r+c-h$. We have
\begin{equation}\label{eq:fan-kernel-threshold}
1-K_q(h,r,c)\le q^d \quad\text{if }d\ge1,
\qquad
K_q(h,r,c)\le q^{1-d} \quad\text{if }d\le0.
\end{equation}
\end{lemma}

\begin{proof}
For (a), at each up-step $g_{j+1}=g_j+1$,
the update $Z_{j+1}=Z_j+1$ occurs with probability
$q^{Z_j-g_j}\le q^{-g_j}$, since $Z_j\ge0$.
For $j\ge M$, this probability is at most $q^{\delta j}$, so a union bound gives an upper bound $q^{\delta M}/(1-q^\delta)$ for the probability of any increase after step $M$. Since this bound tends to zero as $M\to\infty$,
the soft maxima stabilize almost surely.

For (b), we use the identities
\begin{align}
 1-F_h(x)&=xq^{-h}(1-q^h)F_{h-1}(x),
 && h\ge1, \label{eq:F-recurrence-one}\\
 1-(1-x)F_h(x)&=xq^{-h}F_h(qx),
 && h\ge0, \label{eq:F-recurrence-two}
\end{align}
which follow from \eqref{eq:F} by comparing coefficients.

If $d\ge1$, the two identities give
\[
 1-K_q(h,r,c)=
 \begin{cases}
  q^d(1-q^h)F_{h-1}(q^{r+1}),
      &c=1,\ h\ge1,\\
  q^d F_h(q^{r+1}),&c=0,\\
  0,&c=1,\ h=0.
 \end{cases}
\]
Since $0\le F_s(q^a)\le1$ for $a\ge1$ by
Theorem~\ref{thm:longest-kernel}, we obtain
$1-K_q(h,r,c)\le q^d$.

For $d\le0$, applying the first bound to
$(r,h,1-c)$ and using \eqref{eq:kernel-symmetry} gives $K_q(h,r,c)\le q^{1-d}$.

\end{proof}

\subsection{Localization for the auxiliary process} 
We briefly recall the local equilibrium results for ASEP. For general initial conditions, these were established in \cite{AndjelVares1987,BenassiFouque1987}.
See also \cite{Liggett1975,Liggett1977} for earlier results and \cite{Rost1981} for TASEP. In particular, for the step initial condition, if $x_t/t\to v$, then for every fixed $M$,
\begin{equation}\label{eq:fan-step-local-equilibrium}
 \bigl(\eta_t^{x_t}(j)\bigr)_{|j|\le M}
 \Longrightarrow
 \bigl(1-B_j\bigr)_{|j|\le M},
\end{equation}
where the $B_j$ are i.i.d. Bernoulli$(\rho(v))$ random variables.

To use the local equilibrium result in
\eqref{eq:second-class-formula}, we need an additional argument, since local equilibrium gives convergence only on fixed finite intervals, whereas the formula involves soft maxima of infinite paths. The following lemma provides the required localization.

We use the notation for paths from
Section~\ref{sec:second-class}.

\begin{lemma}\label{lem:fan-localization}
Assume \eqref{eq:fan-density-assumption}, let $x_t/t\to v$,
and set $\rho=\rho(v)$. Let $D^\pm$ be defined by
\eqref{eq:full-line-distances} with $x=x_t$.
If $p_R<\rho<p_L$, then
\begin{equation}\label{eq:fan-localization}
 \lim_{M\to\infty}\limsup_{t\to\infty}
 \PP\Bigl(
 (D^-,D^+)\ne
 (D^{-,[M]}(\eta_t^{x_t}),D^{+,[M]}(\eta_t^{x_t}))
 \Bigr)=0.
\end{equation}
The corresponding one-sided conclusion holds
for $D^-$ if $\rho<p_L$ and for $D^+$
if $\rho>p_R$.

\end{lemma}

\begin{proof}
Write $\eta=\eta_t^{x_t}$. For $\delta>0$, let
\[
E_{t,M}=\left\{
\max\{G_n^-(\eta),G_n^+(\eta)\}\le-\delta n
\text{ for all }n\ge M
\right\}.
\]
To prove the lemma, it suffices to show that, for some $\delta>0$,
\begin{equation}\label{eq:fan-negative-tails}
\lim_{M\to\infty}\limsup_{t\to\infty}
\PP(E_{t,M}^c)=0.
\end{equation}
Indeed, on $E_{t,M}$, Lemma~\ref{lem:fan-estimates}(a) and a union bound give
\[
\PP\bigl((D^-,D^+)\ne
(D^{-,[M]}(\eta),D^{+,[M]}(\eta))\bigr)
\le \PP(E_{t,M}^c)+\frac{2q^{\delta M}}{1-q^\delta}.
\]
Thus, \eqref{eq:fan-negative-tails} implies
\eqref{eq:fan-localization}.

To prove \eqref{eq:fan-negative-tails}, let $Y_j(t)$ be the position of color $j$ in the fully colored ASEP started from the identity configuration. By translation invariance of the ASEP, the sequence $W_t(j)=Y_j(t)-j$ is stationary in $j$. By the full-line color-position
symmetry \cite[Theorem~3.1 and Remark~3.3]{BorodinBufetov2021},
\[
(\eta(j))_{j\in\ZZ}
\overset{\mathrm d}=
\bigl(\ind_{\{Y_j(t)\le x_t\}}\bigr)_{j\in\ZZ}.
\]
Since only the distribution of $\eta$ matters, we may realize it as
$\eta(j)=\ind_{\{Y_j(t)\le x_t\}}$ and introduce the stationary sequence $B_t(j)=\ind_{\{W_t(j)>x_t\}}$. Since $1-\eta(j)=\ind_{\{W_t(j)>x_t-j\}}$, for $j\ge1$ we have
\[
1-\eta(-j)\le B_t(-j),
\qquad
1-\eta(j)\ge B_t(j).
\]
Consequently, for every $n\ge1$,
\[
\begin{aligned}
G_n^-(\eta)
&\le\sum_{j=1}^n(B_t(-j)-\zeta(-j)),\\
G_n^+(\eta)
&\le\sum_{j=1}^n(\zeta(j)-B_t(j)).
\end{aligned}
\]
Since $p_R<\rho<p_L$, we can choose $\delta>0$ such that $3\delta<\min\{p_L-\rho,\rho-p_R\}$, so the one-sided density assumptions on $\zeta$ and the bounds for $G_n^\pm$ then reduce \eqref{eq:fan-negative-tails} to proving that, for every $\varepsilon>0$,
\[
\lim_{M\to\infty}\limsup_{t\to\infty}
\PP\left(
\sup_{n\ge M}
\left|\frac1n\sum_{j=1}^n B_t(\pm j)-\rho\right|
>\varepsilon
\right)=0
\]
for either choice of sign. We first consider the positive indices.

We first show that $B_t$ converges locally in distribution to an i.i.d. Bernoulli$(\rho)$ sequence.  By \eqref{eq:fan-step-local-equilibrium}, $1-\eta$ converges locally in distribution to an i.i.d. Bernoulli$(\rho)$ sequence. Thus, to obtain the same convergence for $B_t$, it suffices to show that, for each fixed $j\in\ZZ$,
\[
\PP(B_t(j)\ne 1-\eta(j))\longrightarrow0,
\qquad t\to\infty.
\]
To check this, note that $B_t(j)$ and $1-\eta(j)$ can differ only if $|W_t(j)-x_t|\le |j|$, and using stationarity of $W_t$, color-position symmetry,
and local equilibrium, for every fixed $c\in\ZZ$, we have
\[
\PP(W_t(j)\le x_t+c)
=\PP(\eta_t^{x_t+c}(0)=1)
\longrightarrow1-\rho.
\]
If we take the difference of these probabilities at $c=|j|$ and $c=-|j|-1$, we get
\[
\PP(|W_t(j)-x_t|\le |j|)\longrightarrow0,
\]
as required.

Now we want to control the partial sums
$\frac1n\sum_{j=1}^n B_t(j)$ uniformly for
all sufficiently large $n$. For this purpose, we divide $(B_t(j))_{j\ge1}$ into blocks of fixed length $k\ge1$ and define
\[
Z_{t,i}=
\left|
\frac1k\sum_{j=(i-1)k+1}^{ik}B_t(j)-\rho
\right|,
\qquad i\ge1.
\]
Note that sequence $(Z_{t,i})_{i\ge1}$ is stationary. For any $n\ge M\ge k$, write $n=ak+s$, where $0\le s<k$. The first $ak$ terms form $a$ complete blocks, and the remaining $s$ terms contribute at most $k/M$ to the empirical average. Hence, we have
\[
\left|
\frac1n\sum_{j=1}^nB_t(j)-\rho
\right|
\le
\frac1a\sum_{i=1}^a Z_{t,i}+\frac{k}{M}.
\]
Take $M$ sufficiently large such 
that $k/M<\varepsilon$. Then we have
\[
\begin{aligned}
&\PP\left(
\sup_{n\ge M}
\left|\frac1n\sum_{j=1}^nB_t(j)-\rho\right|
>\varepsilon
\right)\\
&\qquad\le
\PP\left(
\sup_{a\ge1}\frac1a\sum_{i=1}^a Z_{t,i}
>\varepsilon-\frac{k}{M}
\right)\\
&\qquad\le
\frac{\EE Z_{t,1}}{\varepsilon-k/M},
\end{aligned}
\]
where the second inequality follows from the maximal
ergodic inequality
\cite[Chapter~1, Corollary~2.2]{Krengel1985}. So, it remains to estimate $\EE Z_{t,1}$. To do this, fix $k$, and we use the local convergence of $B_t$ and Cauchy--Schwarz inequality
\[
\begin{aligned}
\limsup_{t\to\infty}\EE Z_{t,1}
&\le
\lim_{t\to\infty}
\left[
\EE\left(
\frac1k\sum_{j=1}^k B_t(j)-\rho
\right)^2
\right]^{1/2}=\sqrt{\frac{\rho(1-\rho)}{k}}.
\end{aligned}
\]
Combining the two estimates, we obtain
\[
\begin{aligned}
&\lim_{M\to\infty}\limsup_{t\to\infty}
\PP\left(
\sup_{n\ge M}
\left|\frac1n\sum_{j=1}^n B_t(j)-\rho\right|
>\varepsilon
\right)\le
\frac1\varepsilon
\sqrt{\frac{\rho(1-\rho)}{k}}.
\end{aligned}
\]
Since $k$ is arbitrary, letting $k\to\infty$ shows that the left-hand side is zero. The same argument applies to $(B_t(-j))_{j\ge1}$. This proves \eqref{eq:fan-negative-tails}.
The one-sided conclusions follow by the same argument on the corresponding half-line.

\end{proof}

\subsection{Proof of Theorem~\ref{thm:general-rarefaction}}
We first prove that $\Phi_{q,\zeta}$ is continuous
on $(p_R,p_L)$.
For $M\ge0$, let
$S^{\pm,[M]}(\rho)=(S_j^\pm(\rho))_{j=0}^M$
be the truncated paths and define
\[
\Phi_{q,\zeta}^{[M]}(\rho)
=\EE K_q\bigl(
\soft(S^{-,[M]}(\rho)),
\soft(S^{+,[M]}(\rho)),
1-B_0
\bigr).
\]

Since $\Phi_{q,\zeta}^{[M]}$ is a polynomial in $\rho$,
it suffices to show that, for every compact
$J\subset(p_R,p_L)$,
\begin{equation}\label{eq:fan-static-truncation}
 \sup_{\rho\in J}
 |\Phi_{q,\zeta}(\rho)-\Phi_{q,\zeta}^{[M]}(\rho)|
 \longrightarrow0
 \qquad\text{as }M\to\infty.
\end{equation}
Since $J$ is a compact subset of $(p_R,p_L)$,
\eqref{eq:fan-density-assumption} implies that,
for some $\delta>0$ and all sufficiently large $n$,
\[
 \sup_{\rho\in J}\EE S_n^\pm(\rho)\le-2\delta n.
\]
By Hoeffding's inequality and a union bound,
\[
 \sup_{\rho\in J}
 \PP\bigl(\exists n\ge M:
 S_n^-(\rho)>-\delta n
 \text{ or }S_n^+(\rho)>-\delta n\bigr)
 \longrightarrow0
 \qquad\text{as }M\to\infty.
\]
On the complementary event, Lemma~\ref{lem:fan-estimates}(a),
applied conditionally on the paths, shows that the probability
that either soft maximum changes after step $M$ is at most
$2q^{\delta M}/(1-q^\delta)$.
Since $0\le K_q\le1$, this proves
\eqref{eq:fan-static-truncation}.

Let $x_t/t\to v\in(v_L,v_R)$ and set $\rho=\rho(v)$.
By local equilibrium \eqref{eq:fan-step-local-equilibrium},
for every fixed $M$,
\[
 \EE K_q\bigl(
 D^{-,[M]}(\eta_t^{x_t}),
 D^{+,[M]}(\eta_t^{x_t}),
 \eta_t^{x_t}(0)\bigr)
 \longrightarrow\Phi_{q,\zeta}^{[M]}(\rho)
 \qquad\text{as }t\to\infty.
\]
By \eqref{eq:second-class-formula} and $0\le K_q\le1$,
replacing $D^\pm$ by their finite-prefix versions
changes the expectation by at most
\[
\PP\Bigl(
(D^-,D^+)\ne
(D^{-,[M]}(\eta_t^{x_t}),
 D^{+,[M]}(\eta_t^{x_t}))
\Bigr).
\]
By \eqref{eq:fan-localization}, the limsup of this error
as $t\to\infty$ tends to zero as $M\to\infty$.
Together with the convergence above and
\eqref{eq:fan-static-truncation}, this gives
\begin{equation}\label{eq:fan-interior-convergence}
 \PP(Q_\zeta(t)\le x_t)
 \longrightarrow\Phi_{q,\zeta}(\rho(v)).
\end{equation}

It remains to consider the behavior outside the fan.
Suppose $v<v_L$, set $x_t=\lfloor vt\rfloor$,
and write $\rho=\rho(v)$.
If $\rho>p_L$, the strong law of large numbers and
\eqref{eq:fan-density-assumption} give, almost surely,
\[
\frac{S_n^-(\rho)}{n}\longrightarrow \rho-p_L>0.
\]
By \eqref{eq:running-record}, it follows that
$\soft(S^{-,[M]}(\rho))\to\infty$ almost surely
as $M\to\infty$.
For each fixed $M$, local equilibrium gives
\[
D^{-,[M]}(\eta_t^{x_t})
\Longrightarrow \soft(S^{-,[M]}(\rho))
\qquad\text{as }t\to\infty.
\]
Since $D^-\ge D^{-,[M]}(\eta_t^{x_t})$,
we obtain $D^-\to\infty$ in probability as $t\to\infty$.
On the other hand, since $\rho>p_R$,
local equilibrium and the one-sided version of
Lemma~\ref{lem:fan-localization} imply that
$D^+$ is tight as $t\to\infty$.
Thus, by \eqref{eq:second-class-formula} and
\eqref{eq:fan-kernel-threshold},
\[
\PP(Q_\zeta(t)\le x_t)\longrightarrow0.
\]
If $\rho=p_L=1$, the path $S^+(1)$ has no up-steps.
Hence, local equilibrium and one-sided localization give
$D^+\to0$ and $\eta_t^{x_t}(0)\to0$
in probability as $t\to\infty$.
Since $K_q(h,0,0)=0$, we again obtain
$\PP(Q_\zeta(t)\le x_t)\to0$ as $t\to\infty$.

The case $v>v_R$ is analogous. Set $x_t=\lfloor vt\rfloor$ and $\rho=\rho(v)$.
If $\rho<p_R$, the strong law of large numbers and \eqref{eq:fan-density-assumption} give, almost surely,
\[
\frac{S_n^+(\rho)}{n}\longrightarrow p_R-\rho>0
\qquad\text{as }.
\]
As above, local equilibrium and prefix monotonicity imply that $D^+\to\infty$ in probability as $t\to\infty$, while one-sided localization gives tightness of $D^-$. By \eqref{eq:second-class-formula} and
\eqref{eq:fan-kernel-threshold},
\[
\PP(Q_\zeta(t)\le x_t)\longrightarrow1.
\]
If $\rho=p_R=0$, the path $S^-(0)$ has no up-steps.
Hence, local equilibrium and one-sided localization give
$D^-\to0$ and $\eta_t^{x_t}(0)\to1$
in probability as $t\to\infty$.
Since $K_q(0,r,1)=1$, we obtain the same limit.

Together with \eqref{eq:fan-interior-convergence}
and the continuity of $\Phi_{q,\zeta}$,
this proves \eqref{eq:fan-weak-limit}.

\section{Explicit examples}\label{sec:rarefaction-examples}
We now apply Theorem~\ref{thm:general-rarefaction}
to derive explicit limiting distributions of the
second-class particle for different classes of initial conditions.

\subsection{Technical lemmas}

We begin with two technical lemmas used in the explicit calculations below.

In the original construction of the soft maximum \eqref{eq:soft-scan}, the update probabilities depend on both the current state and the height of the path. Lemma~\ref{lem:reverse-queue} provides a reverse-order sampling rule that removes the dependence on the path height.

Lemma~\ref{lem:iid-soft-law} applies this rule to paths with i.i.d. increments, which arise from certain classes of initial conditions, and gives an explicit distribution for their soft maxima.

\begin{lemma}\label{lem:reverse-queue}
Let $g$ be a finite discrepancy path with increments
$e_1,\ldots,e_n$.
Starting from $N=0$, read the increments in reverse order
$e_n,\ldots,e_1$ and update $N$ as follows.
For an increment $+1$, increase $N$ by one.
For an increment $-1$, decrease $N$ by one with
probability $1-q^N$, leaving it unchanged otherwise.
For an increment $0$, leave $N$ unchanged.
All random choices are independent.

The final value of $N$ has the same distribution
as $\soft(g)$.
\end{lemma}

\begin{proof}
The polynomials \eqref{eq:packed-polynomial} satisfy
\[
\poly_{N+1}(u)=(1+u)\poly_N(u/q)
\]
and, for $N\ge1$,
\[
q^N\poly_N(u)+(1-q^N)\poly_{N-1}(u)
=\poly_N(qu).
\]
Let $N_g$ denote the output of the reverse algorithm.
By Proposition~\ref{prop:transform}, it suffices to show that
$\EE\poly_{N_g}(u)=\Pi_g(u)$, which we prove by induction on the length of $g$. If $g$ is the empty path, the claim is clear.
Now consider a path $g$ of positive length, and let
$\widetilde g$ be the path starting at zero with increments
$e_2,\ldots,e_n$. 
Since the reverse algorithm reads $e_1$ last,
there are three cases to consider. If $e_1=-1$, then conditionally on
$N_{\widetilde g}=k$, for $k\ge1$, using the second
identity above, we obtain
\[
\begin{aligned}
\EE[\poly_{N_g}(u)\mid N_{\widetilde g}=k]
&=q^k\poly_k(u)+(1-q^k)\poly_{k-1}(u)\\
&=\poly_k(qu).
\end{aligned}
\]
Taking expectations and using the induction hypothesis,
we obtain
$\EE\poly_{N_g}(u)=\Pi_{\widetilde g}(qu)$.
The cases $e_1=0$ and $e_1=1$ are checked similarly,
but for $e_1=1$ we use the other polynomial identity. Thus, we obtain
\[
\EE\poly_{N_g}(u)=
\begin{cases}
(1+u)\Pi_{\widetilde g}(u/q),&e_1=1,\\
\Pi_{\widetilde g}(u),&e_1=0,\\
\Pi_{\widetilde g}(qu),&e_1=-1.
\end{cases}
\]
The right-hand side equals $\Pi_g(u)$ in each case,
which completes the proof.
\end{proof}

For $0\le a<1$, let $\pi_a$ denote the
$q$-geometric distribution given by
\begin{equation}\label{eq:q-geometric}
\pi_a(n)=(a;q)_\infty\frac{a^n}{(q;q)_n},
\qquad n\in\NN,
\qquad (a;q)_\infty=\prod_{j\ge0}(1-aq^j).
\end{equation}

\begin{lemma}\label{lem:iid-soft-law}
Let $g$ be an infinite path starting at zero with
i.i.d. increments in $\{-1,0,1\}$.
Suppose that the probabilities of increments $+1$
and $-1$ are $u$ and $d$, respectively,
where $0\le u<d$.
Then $\soft(g)$ has distribution $\pi_{u/d}$.

Moreover, if $H$ and $R$ are independent random variables
with distributions $\pi_a$ and $\pi_b$, respectively,
where $0\le a,b<1$, then
\begin{equation}\label{eq:q-geometric-kernel}
\EE K_q(H,R,1)=\frac{1-a}{1-ab},\qquad
\EE K_q(H,R,0)=\frac{b(1-a)}{1-ab}.
\end{equation}
\end{lemma}

\begin{proof}
By the strong law of large numbers, $g_n/n\to u-d<0$
almost surely, so the finite-prefix soft maxima of $g$
stabilize almost surely by Lemma~\ref{lem:fan-estimates}(a).

For each $m\ge0$, by Lemma~\ref{lem:reverse-queue},
the soft maximum of the first $m$ increments of $g$
has the same distribution as a Markov chain after $m$ steps,
started at zero, with transition probabilities
$P(n,n+1)=u$ and $P(n,n-1)=d(1-q^n)$.

For this chain, the distribution $\pi_{u/d}$ satisfies
the detailed balance equations
\[
(1-q^n)\pi_a(n)=a\pi_a(n-1),\qquad n\ge1.
\]
For $u>0$, our Markov chain is irreducible and aperiodic,
so $\pi_{u/d}$ is its unique stationary distribution
and the chain converges to it,
while for $u=0$, the chain remains at zero.
Thus, $\soft(g)$ has distribution $\pi_{u/d}$.

For the second part of the lemma, let $J(a,b)=\EE K_q(H,R,1)$.
By \eqref{eq:kernel},
$K_q(h,r,0)=(1-q^r)K_q(h,r-1,1)$ for $r\ge1$,
while $K_q(h,0,0)=0$.
Using the detailed balance identity
$(1-q^r)\pi_b(r)=b\pi_b(r-1)$, we obtain
\[
\begin{aligned}
\EE K_q(h,R,0)
&=\sum_{r\ge1}\pi_b(r)(1-q^r)K_q(h,r-1,1)\\
&=b\sum_{r\ge0}\pi_b(r)K_q(h,r,1)\\
&=b\,\EE K_q(h,R,1).
\end{aligned}
\]
Since $H$ and $R$ are independent, averaging over $H$ gives $\EE K_q(H,R,0)=bJ(a,b)$. To find $J(a,b)$, we use the kernel symmetry \eqref{eq:kernel-symmetry}, which gives
\[
\begin{aligned}
J(a,b)
&=1-\EE K_q(R,H,0)\\
&=1-aJ(b,a).
\end{aligned}
\]
Interchanging $a$ and $b$, we obtain
$J(b,a)=1-bJ(a,b)$, and hence
\[
J(a,b)=\frac{1-a}{1-ab}.
\]
\end{proof}

\subsection{Bernoulli distributions on both half-lines}
\begin{proposition}\label{prop:bernoulli-fan}
Let $\zeta(y)$, $y\ne0$, be independent Bernoulli random variables with parameters $p_1$ for $y<0$
and $p_2$ for $y>0$, where $0\le p_2<p_1\le1$.
Then
\begin{equation}\label{eq:bernoulli-uniform}
\frac{Q_\zeta(t)}t\Longrightarrow
\operatorname{Unif}\bigl[
\gamma(1-2p_1),\gamma(1-2p_2)\bigr].
\end{equation}
\end{proposition}

\begin{proof}
For $p_2<\rho<p_1$, the increments of the two paths
in \eqref{eq:fan-static-paths} are i.i.d., with
upward and downward probabilities
$\rho(1-p_1),(1-\rho)p_1$ on the left and
$(1-\rho)p_2,\rho(1-p_2)$ on the right.
By Lemma~\ref{lem:iid-soft-law}, the two soft maxima
are independent with distributions $\pi_a$ and $\pi_b$,
where
\begin{equation}\label{eq:bernoulli-ab}
a=\frac{\rho(1-p_1)}{(1-\rho)p_1},\qquad
b=\frac{(1-\rho)p_2}{\rho(1-p_2)}.
\end{equation}
Since $B_0$ is independent of the two soft maxima,
\eqref{eq:q-geometric-kernel} gives
\begin{equation}\label{eq:bernoulli-functional}
\EE_\zeta\Phi_{q,\zeta}(\rho)
=(1-\rho)\frac{1-a}{1-ab}
+\rho\frac{b(1-a)}{1-ab}
=\frac{p_1-\rho}{p_1-p_2}.
\end{equation}
By Theorem~\ref{thm:general-rarefaction},
the limiting distribution is uniform.
\end{proof}
This produces an (admittedly, very complicated) independent proof of the results of Ferrari-Kipnis
\cite{FerrariKipnis1995} and Ferrari-Goncalves-Martin
\cite{FerrariGoncalvesMartin2009}. 

\subsection{Iid gaps on the left half-axis}
In this section we consider a fairly general class of initial configurations on the left half-axis, while the right half-axis will be assumed to be empty. 
\begin{proposition}\label{prop:gaps-fan}
Let $A,G_1,G_2,\ldots$ be independent positive integer-valued random variables. Assume that $G_1,G_2,\ldots$ are identically distributed, and denote a random variable with their common distribution by $G$. Consider the initial condition with first-class particles at
\begin{equation}\label{eq:gaps-initial}
-A,\quad -A-G_1,\quad -A-G_1-G_2,\quad\ldots,
\end{equation}
a second-class particle at zero, and holes elsewhere. Set $\mu=\EE G<\infty$ and define the probability generating functions
\begin{equation}\label{eq:gaps-generating}
\psi_A(s)=\EE s^{A-1},\qquad
\psi_G(s)=\EE s^{G-1}.
\end{equation}
Then the limiting distribution of the second-class
particle is supported on
$[\gamma(1-2/\mu),\gamma]$, and
\begin{equation}\label{eq:gaps-fan}
\lim_{t\to\infty}
\PP\left(\frac{Q_\zeta(t)}{t}\le v\right)
=(1-\mu\rho(v))
\frac{\psi_A(1-\gamma\rho(v))}
{\psi_G(1-\gamma\rho(v))},
\qquad \gamma(1-2/\mu)<v<\gamma.
\end{equation}
\end{proposition}
\begin{proof}
By the strong law of large numbers, $p_L=1/\mu$ and $p_R=0$ almost surely. 
Fix $0<\rho<1/\mu$, so the strong law of large numbers and Lemma~\ref{lem:fan-estimates}(a) imply that the finite-prefix soft maxima of the left path stabilize almost surely. The right path has no up-steps, so $R_\rho=0$. Since $K_q(h,0,1)=q^h$ and $K_q(h,0,0)=0$, we have \begin{equation}\label{eq:empty-right-functional}
    \EE_\zeta\Phi_{q,\zeta}(\rho)
    =(1-\rho)\EE q^{H_\rho}.
\end{equation}
Thus, it remains to compute the $q$-moment
$\EE q^{H_\rho}$.

We first consider the case $A=1$.
The first-class particles are located at
\[
-1,\quad -1-G_1,\quad -1-G_1-G_2,\quad\ldots.
\]
If we delete the $G_1$ sites nearest the origin and shift all remaining sites by $G_1$ steps to the right, the first-class particles are then located at
\[
-1,\quad -1-G_2,\quad -1-G_2-G_3,\quad\ldots.
\]
Since $G_1,G_2,\ldots$ are i.i.d., the resulting configuration has the same law as the original one and is independent of $G_1$.
 
For this new configuration, run the reverse-order algorithm of Lemma~\ref{lem:reverse-queue}, and denote its output by $N$. Since the two configurations have the same law, $N\overset{\mathrm d}=H_\rho$, and $N$ is also independent of $G_1$.

Now consider the reverse-order algorithm for the original configuration. After reading the remaining configuration, the current state is $N$. So, it remains to read the $G_1-1$ holes and the particle at $-1$. Each of the $G_1-1$ holes gives an increment $+1$ in the left discrepancy path \eqref{eq:fan-static-paths} with probability $\rho$, and $0$ otherwise, while the particle at $-1$ gives an increment $-1$ with probability $1-\rho$, and $0$ otherwise. In the reverse-order algorithm, the holes increase the state from $N$ to $N+L$, where the variable $L$ has distribution $\operatorname{Bin}(G_1-1,\rho)$ and is independent of $N$. 
The particle at $-1$ then decreases the state by $D\in\{0,1\}$, where $D$ is Bernoulli with parameter $(1-\rho)(1-q^{N+L})$. Thus, the output of the algorithm for the original
configuration is
\begin{equation}\label{eq:gaps-reverse-cycle}
    N'=N+L-D\overset{\mathrm d}=N,
\end{equation} 
since both $N$ and $N'$ have the same distribution as $H_\rho$.

We now show that $\EE L=\EE D$. Since $\EE N$ may be infinite, we use a standard
truncation argument. By \eqref{eq:gaps-reverse-cycle}, for every $K\ge1$,
$
\EE[(N'\wedge K)-(N\wedge K)]=0$. Moreover, $|N'-N|\le L+1\le G_1$, so dominated convergence gives $\EE(N'-N)=0$, and hence $\EE L=\EE D$.

Since $L$ has distribution
$\operatorname{Bin}(G_1-1,\rho)$ conditionally on $G_1$,
we have $\EE L=\rho(\mu-1)$ and
$\EE q^L =\EE(1-\rho+\rho q)^{G_1-1}$, which equals to $\psi_G(1-\gamma\rho)$.
Since $D$ is Bernoulli with parameter
$(1-\rho)(1-q^{N+L})$ conditionally on $N$ and $L$, we have
\[
\begin{aligned}
\EE D
&=(1-\rho)\EE(1-q^{N+L})\\
&=(1-\rho)(1-\psi_G(1-\gamma\rho)\EE q^N).
\end{aligned}
\]
Thus, $\EE L=\EE D$ gives
\begin{equation}\label{eq:gaps-cycle-balance}
\EE q^N
=\frac{1-\mu\rho}{(1-\rho)\psi_G(1-\gamma\rho)},
\end{equation}
which equals $\EE q^{H_\rho}$ for $A=1$.

For $A >1$, the first-class particles are located at
\[
-A,\quad -A-G_1,\quad -A-G_1-G_2,\quad\ldots.
\]
In the reverse-order algorithm, reading the configuration up to site $-A$ gives a state with the same distribution as $N$, since this part of the configuration is a shift
of the case $A=1$. Thus, it remains to read the $A-1$ holes at sites $-A+1,\ldots,-1$. 
These holes contribute an increment $L_A$,
independent of $N$, which has distribution
$\operatorname{Bin}(A-1,\rho)$. Thus, the output of the algorithm has the same distribution as $N+L_A$, and hence
\[
\EE q^{H_\rho}
=\EE q^N\,\EE q^{L_A}
=\EE q^N\,\psi_A(1-\gamma\rho).
\]
Combining this with \eqref{eq:empty-right-functional} and \eqref{eq:gaps-cycle-balance}, we obtain
\[
\EE_\zeta\Phi_{q,\zeta}(\rho)
=(1-\mu\rho)
\frac{\psi_A(1-\gamma\rho)}
{\psi_G(1-\gamma\rho)}.
\]
\end{proof}
If $A$ and $G$ have the same distribution, then
$\psi_A=\psi_G$, and \eqref{eq:gaps-fan} gives
\[
1-\mu\rho(v)
=\frac{v-\gamma(1-2/\mu)}{2\gamma/\mu}.
\]
Thus, the limiting distribution is uniform on the fan. In particular, if $A$ and $G$ have the geometric distribution
\[
    \PP(A=k)=\PP(G=k)=p(1-p)^{k-1},
    \qquad k\ge1,
\]
we recover the Bernoulli initial condition with density $p$ on the left and an empty right half-line.
\begin{example}\label{ex:gaps-one-three}
Let $A=1$ and suppose that
$\PP(G=1)=\PP(G=3)=1/2$.
Then $\mu=2$, $\psi_A(s)=1$, and
$\psi_G(s)=(1+s^2)/2$.
The limiting distribution is supported on $[0,\gamma]$
and has distribution function
\begin{equation}\label{eq:gaps-one-three}
\begin{aligned}
\lim_{t\to\infty}\PP\left(\frac{Q_\zeta(t)}t\le v\right)
&=\frac{2(1-2\rho(v))}{1+[1-\gamma\rho(v)]^2}\\
&=\frac{2v/\gamma}{1+\bigl((1+q+v)/2\bigr)^2},
\qquad 0<v<\gamma.
\end{aligned}
\end{equation}
\end{example}
An important special case of Proposition~\ref{prop:gaps-fan} is given by particles arranged in an arithmetic progression.
\begin{example}\label{ex:periodic-fan}
Fix integers $m,d\ge1$. Consider the initial condition with first-class particles at
\[
-d,-d-m,-d-2m,\ldots,
\]
a second-class particle at zero, and holes elsewhere. Taking $A=d$ and $G=m$ in
Proposition~\ref{prop:gaps-fan}, we obtain the limiting distribution supported on
$[\gamma(1-2/m),\gamma]$, with distribution function
\begin{equation}\label{eq:periodic-fan}
F_{m,d}(v)
=(1-m\rho(v))[1-\gamma\rho(v)]^{d-m},
\qquad \gamma(1-2/m)<v<\gamma.
\end{equation}
\end{example}
As usual in this type of problems, different limiting distributions can arise even with the same particle density $1/m$ on the left. For instance, particles at $-2,-4,-6,\ldots$ give the uniform distribution on $[0,\gamma]$, whereas particles at $-1,-3,-5,\ldots$ give
\begin{equation}\label{eq:alternating-phase}
F_{2,1}(v)
=\frac{1-2\rho(v)}{1-\gamma\rho(v)},
\qquad 0<v<\gamma.
\end{equation}

\subsection{A one-site perturbation of a Bernoulli background}
We now show that even changing the occupation probability at a single site in a Bernoulli initial condition can affect the limiting distribution of the
second-class particle.
\begin{proposition}\label{prop:bernoulli-defect}
Let $0\le p_2<p_1\le1$ and $\theta\in[0,1]$.
Let $\zeta(y)$, $y\ne0$, be independent Bernoulli random variables with parameters $p_1$ for $y\le-2$, $\theta$ for $y=-1$, and $p_2$ for $y\ge1$. For $p_2<\rho(v)<p_1$, set
\begin{equation}\label{eq:defect-coordinate}
    w(v)=\frac{p_1-\rho(v)}{p_1-p_2}.
\end{equation}
The limiting distribution is supported on
$[\gamma(1-2p_1),\gamma(1-2p_2)]$ and has
distribution function
\begin{equation}\label{eq:defect-cdf}
    F_\theta(v)
    =w(v)+
    \frac{\gamma(\theta-p_1)(p_1-p_2)}
    {p_1(1-p_2)-qp_2(1-p_1)}
    w(v)(1-w(v))
\end{equation}
in the interior of the fan.
\end{proposition}
\begin{proof}
Write $\rho=\rho(v)$ and denote by $F_0(v)$ and $F_1(v)$ the limiting distribution functions when site $-1$ is empty and occupied, respectively.
Since the occupation at site $-1$ is independent of the remaining configuration, we have
\[
    F_\theta(v)=(1-\theta)F_0(v)+\theta F_1(v).
\]
The case $\theta=p_1$ corresponds to Bernoulli initial conditions on both half-lines, for which we have already shown that
\[
    w(v)=(1-p_1)F_0(v)+p_1F_1(v).
\]
Combining these two identities, we get
\[
    F_\theta(v)
    =F_0(v)+\frac{\theta}{p_1}
    \bigl(w(v)-F_0(v)\bigr).
\]
Thus it suffices to compute $F_0(v)$.

Let $H$ and $R$ be the soft maxima of the left path starting at site $-2$ and the right path, respectively. Since site $-1$ is empty, its contribution to the left path \eqref{eq:fan-static-paths} is the increment
$B_{-1}$. In the reverse-order algorithm, this increment is processed last, and by Lemma~\ref{lem:reverse-queue}, we obtain $H_\rho\overset{\mathrm d}=H+B_{-1}$. Thus, by Theorem~\ref{thm:general-rarefaction},
\[
    F_0(v)=\EE K_q(H+B_{-1},R, 1-B_0).
\]
The pair $(H,R)$ has the same distribution as the corresponding soft maxima in the Bernoulli case. In particular, $H$ and $R$ have distributions $\pi_a$ and $\pi_b$, respectively, where $a,b$ are given by \eqref{eq:bernoulli-ab}. Moreover, the computation in the proof of Lemma~\ref{lem:iid-soft-law} gives $\EE_R K_q(h,R,0)=b\,\EE_R K_q(h,R,1)$ for every $h$. Using this identity and averaging over $B_0$, we obtain
\[
    F_0(v)=(1-\rho+\rho b)
    \EE K_q(H+B_{-1},R,1).
\]
Next, averaging over $B_{-1}$, we obtain
\begin{equation}\label{eq:defect-hole-functional}
    F_0(v)=(1-\rho+\rho b)\bigl[
    (1-\rho)\EE K_q(H,R,1)
    +\rho\EE K_q(H+1,R,1)\bigr].
\end{equation}
The first expectation is known from
Lemma~\ref{lem:iid-soft-law}, namely
$\EE K_q(H,R,1)=(1-a)/(1-ab)$. Thus it remains to compute $\EE K_q(H+1,R,1)$. 

For $a>0$, we have
\[
\begin{aligned}
    \EE K_q(H+1,R,1)
    &=\sum_{n\ge1}\pi_a(n-1)\EE_R K_q(n,R,1)\\
    &=\frac1a\sum_{n\ge0}
    \bigl(\pi_a(n)-(1-a)\pi_{qa}(n)\bigr)
    \EE_R K_q(n,R,1)\\
    &=\frac{1-a}{a}
    \left(\frac1{1-ab}-\frac{1-qa}{1-qab}\right).
\end{aligned}
\]
Here second equality uses $\pi_a(n-1)=(1-q^n)\pi_a(n)/a$ and $q^n\pi_a(n)=(1-a)\pi_{qa}(n)$. For the last equality, we split the sum into two and apply Lemma~\ref{lem:iid-soft-law} with
parameters $(a,b)$ and $(qa,b)$, respectively. Simplifying the expression, we get
\begin{equation}\label{eq:defect-shifted-kernel}
    \EE K_q(H+1,R,1)
    =\frac{1-a}{1-ab}
    -\frac{\gamma(1-a)(1-b)}
    {(1-ab)(1-qab)}.
\end{equation}
 Substituting this into \eqref{eq:defect-hole-functional}, using \eqref{eq:bernoulli-functional}, and simplifying, we obtain
\[
    F_0(v)=w(v)-
    \frac{\gamma p_1(p_1-p_2)}
    {p_1(1-p_2)-qp_2(1-p_1)}
    w(v)(1-w(v)).
\]
Combining this with the expression for $F_\theta(v)$ obtained above gives \eqref{eq:defect-cdf} for $a>0$.

The case $a=0$ follows by continuity, since
$\pi_a$ converges to $\delta_0$ as $a\downarrow0$, $K_q$ is bounded, and the right-hand side of \eqref{eq:defect-shifted-kernel} is continuous at $a=0$.
\end{proof}
\Needspace{8\baselineskip}

\end{document}